\documentclass[12pt]{article}
\usepackage{amssymb,amsmath,color}
\def \R{{\hbox{\vrule width 0.6pt height 6.8pt depth -.2pt\kern-0.2pt
R}}}
\def \P{{\hbox{\vrule width 0.6pt height 6.8pt depth -.2pt\kern-0.2pt
P}}}

\newcommand{\aref}[1]{(\ref{#1})}

\def \no { \noindent}
\def \er {\mathbb R}

\def \m {\eta}

\def\grad{\nabla}
\def\div{\, \mbox{div}\,  }

\def\MOD#1{{|\kern -.16em |\kern -.16em | #1 | \kern -.16em |\kern
 -.16em |}}
\def \epsilon {\varepsilon}

\def\ds{\displaystyle}
\newtheorem{theo}{\bf THEOREM}[section]
\newtheorem{lem}[theo]{\bf LEMMA}
\newtheorem{pro}[theo]{\bf PROPOSITION}
\newtheorem{cor}[theo]{\bf COROLLARY}

\numberwithin{equation}{section}

\def\B1{B_{1/2}}
\def\Box{\hfill\rule{2.5mm}{2.5mm}}

\def\R{{\mathbb {R}}}

\def\build#1_#2^#3{\mathrel{
\mathop{\kern 0pt#1}\limits_{#2}^{#3}}}

\def\h1{\mathop{\rm H^1_{\rm loc,\rm u}}}

\def\l2{\mathop{\rm L^2_{\rm loc,\rm u}}}

\def\t0{{t}}

\def\Box{\hfill\rule{2.5mm}{2.5mm}}
\begin{document}
\bibliographystyle{siam}

\title{
A Lyapunov functional and  blow-up results for a class of perturbations  for   semilinear wave
equations in the  critical case }
\author{M. Hamza\\
{\it \small Facult\'e des Sciences de Tunis}\\
H. Zaag\footnote{This author is supported by a grant from the french Agence Nationale de la Recherche, project
ONDENONLIN, reference ANR-06-BLAN-0185.} \\
{\it \small CNRS LAGA Universit\'e Paris 13}}

\maketitle
\begin{abstract}
\no  We  consider in this paper some class of perturbation for the semilinear wave equation with critical
 (in the conformal transform sense) power nonlinearity.
  Working in the framework of  similarity variables, we find a Lyapunov functional for the problem. Using  a two-step argument based on interpolation and a critical Gagliardo-Nirenberg inequality, we show that the blow-up rate of any sigular solution is given by the solution of the non perturbed associated ODE, namely $u''=u^p$.
\end{abstract}

\vspace{0.5cm}


 \vspace{0.4cm}

\maketitle \maketitle
\section{\bf Introduction}
This paper is devoted to the study of blow-up solutions for the following semilinear  wave
equation:
\begin{equation}
\label{1} \left\lbrace
\begin{array}{l}
 \partial_{tt}u=
\Delta u+ |u|^{p-1} u+f(u)+g(\partial_t u),\qquad  (x,t)\in \er^N\times [0,T), \\
\\
(u(x,0),\partial_tu(x,0))=(u_0(x),u_1(x))\in   H^{1}_{loc,u}(\er^N)\times L^{2}_{loc,u}(\er^N),\\
\end{array}
\right.
\end{equation}
with  critical
 power nonlinearity
\begin{equation}\label{pc}
p=p_c\equiv 1+\frac4{N-1}, \quad {\textrm{where}} \quad N\ge 2.
\end{equation}
We assume that the functions
$f$ and  $g$  are locally Lipschitz-continuous  satisfying
the following conditions
\begin{eqnarray*}
(H_f)&& |{f(x)}|\le M(1+|x|^q) \qquad{{\textrm {with}}}\ \ (q<p,\ \ M>0),\\
\\
(H_g)&& |{g(x)}|\le M(1+|x|).\qquad\qquad
\end{eqnarray*}
 The spaces  $L^{2}_{loc,u}(\er^N)$ and $H^{1}_{loc,u}(\er^N)$ are  defined by
\begin{equation*}
L^{2}_{loc,u}(\er^N)=\{u:\er^N\mapsto \er/ \sup_{a\in \er^N}(\int_{|x-a|\le 1}|u(x)|^2dx)<+\infty \},
\end{equation*}
and
\begin{equation*}
H^{1}_{loc,u}(\er^N)=\{u\in L^{2}_{loc,u}(\er^N),|\grad u|\in L^{2}_{loc,u}(\er^N) \}.
\end{equation*}

\bigskip
\no
The  Cauchy problem  of  equation ({\ref{1}}) is  solved
in $H^{1}_{loc,u}\times L^{2}_{loc,u}$. This follows from the finite
speed of propagation and the the wellposdness in $H^{1} \times L^{2}$, valid whenever
$ 1< p<1+\frac{4}{N-2}$ .
The existence of blow-up   solutions for  the  associated  ordinary differential equation  of
(\ref{1})  is a classical result. By using  the
finite
speed of propagation, we conclude that there exists a
 blow-up solution $u(t)$ of (\ref{1}) which depends non trivially on the space variable.
In this paper, we consider  a blow-up solution $u(t)$ of (\ref{1}), we define (see for example
Alinhac \cite{A} and \cite{A1})  $\Gamma$ as the graph of a function $x \mapsto T(x)$ such
that the domain of definition of  $u$  is given by 
$$D_u=\{(x,t)\ \ \big |t<T(x)\}.$$
The set $D_u$ is called the maximal influence domain of $u$. Moreover, from the finite speed of
propagation, $T$ is a $1$-Lipschitz function.
Let $\overline{T}$ be the minimum of $T(x)$ for all $x\in \er^N$. The time
$\overline{T}$  and the graph $\Gamma$  are called (respectively) the blow-up
time and the blow-up graph of $u$. \\
Let us first  introduce the following non-degeneracy condition for $\Gamma$. If we introduce
for all $x\in\er^N$, $t\le T(x)$ and $\delta>0$, the cone
\begin{equation} \label{cone}
C_{x,t,\delta}=\{(\xi,\tau)\neq(x,t)|0\le \tau\le t-\delta|\xi-x|\},
\end{equation}
then our non degeneracy condition is the following: $x_0$ is a non characteristic point if
\begin{equation} \label{cone1}
\exists\delta_0=\delta_0(x_0)\in (0,1) \ \textrm{ such that}\  u \ \textrm{ is defined on }\
C_{x_0,T(x_0),\delta_0}.
\end{equation}
\no
We aim at studying  the growth estimate  of $u(t)$ near the space-time blow-up graph in the critical
case (where $p=p_c$ satisfies (\ref{pc})).

\bigskip
\no
 In the case $(f,g)\equiv (0,0)$, equation (\ref{1}) reduces to the semilinear wave equation:
\begin{equation}\label{mu=0}
\partial_{tt}u=\Delta u+ |u|^{p-1} u,\qquad  (x,t)\in \er^N\times [0,\overline{T}).
\end{equation}

\no  Merle and Zaag   in \cite{MZ3} (see also    \cite{MZ1} and \cite{MZ2}) have proved, that  if
$1<p\le p_c=1+\frac4{N-1}$, if $u$ is a solution of (\ref{mu=0})
with blow up graph
$\Gamma:\{x\mapsto T(x)\}$, then for all $x_0 \in \er^N$  and $t\in [\frac34T(x_0),T(x_0)]$,
  the growth estimate near the space-time blow-up graph  satisfies
\begin{eqnarray*}
 (T(x_0)-t)^{\frac{2}{p-1}}\frac{\|u(t)\|_{L^2(B(x_0,\frac{T(x_0)-t}{2}))}}{ (T(x_0)-t)^{\frac{N}{2}}}\qquad
 \qquad\qquad\qquad\qquad\qquad \qquad\qquad\qquad\\
 + (T(x_0)-t)^{\frac{2}{p-1}+1}\!\Big (\frac{\|\partial_t u(t)\|_{L^2(B(x_0,\frac{T(x_0)-t}{2}))}}{ (T(x_0)-t)^{\frac{N}{2}}}+
 \frac{\|\grad u(t)\|_{L^2(B(x_0,\frac{T(x_0)-t}{2}))}}{ (T(x_0)-t)^{\frac{N}{2}}}\Big )\le K,
\end{eqnarray*}
where the constant $K$  depends only on $N,p,$ and on  an upper bound on $T(x_0)$, $\frac1{T(x_0)}$
and the initial data in $ H^{1}_{loc,u}(\er^N)\times L^{2}_{loc,u}(\er^N)$.
If in addition $x_0$ is non characteristic (in the sense (\ref{cone1})), then for all $t\in [ \frac{3T(x_0)}4, T(x_0)]$,
\begin{eqnarray}\label{mzmz}
0<\varepsilon_0(N,p)\le (T(x_0)-t)^{\frac{2}{p-1}}\frac{\|u(t)\|_{L^2(B(x_0,{T(x_0)-t}))}}{ (T(x_0)-t)^{\frac{N}{2}}}
\qquad \qquad\qquad\qquad\qquad
\nonumber\\
 + (T(x_0)-t)^{\frac{2}{p-1}+1}\Big (\frac{\|\partial_t u(t)\|_{L^2(B(x_0,{T(x_0)-t}))}}{ (T(x_0)-t)^{\frac{N}{2}}}+
 \frac{\|\grad u(t)\|_{L^2(B(x_0,{T(x_0)-t}))}}{ (T(x_0)-t)^{\frac{N}{2}}}\Big )\le K,
\end{eqnarray}
where the constant $K$  depends only on $N,p,$ and on  an upper bound on $T(x_0)$, $\frac1{T(x_0)}$, $\delta_0(x_0)$
and the initial data in $ H^{1}_{loc,u}(\er^N)\times L^{2}_{loc,u}(\er^N)$.

\bigskip

\no Following this blow-up rate estimate, Merle and Zaag addressed the question of the asymptotic
behavior of $u(x,t)$ near $\Gamma$
in one space dimension.

\no More precisely, they proved in \cite{MZ4}
 and \cite{MZ5} that the set of non charecteristic points ${\cal R}\subset \er$ is non empty open
and that $x\mapsto T(x)$ is of class $C^1$ on ${\cal R}$. They also described the blow-up profile of $u$  near $(x_0,T(x_0))$ when $x_0\in {\cal R}$.

\no In \cite{MZarxiv}, they proved that $S=\er\backslash {\cal R}$ has an empty interior and that $\Gamma$ is a corner  of angle $\frac{\pi}{2}$ near any $x_0\in S$. They also showed that
$u(x,t)$ decomposes in a sum  of decoupled solitons near $(x_0,T(x_0))$.
They also gave  examples of blow-up solutions with $S \neq\emptyset$.

\bigskip
\no  In \cite{HZ}, we addressed the question
of extending  the results of Merle and Zaag
\cite {MZ1}, \cite{MZ2} and  \cite{MZ3} to perturbed equations of the type (\ref{1}).
In  \cite{HZ}, we could prove the statement
(\ref{mzmz}) under some reasonable growth estimates an $f$ and $g$ in (\ref{1}) however only when $1<p<p_c$. When $N\ge 2$ and $p=p_c$ our method in \cite{HZ} breaks down. Let us briefly explain in the following how the method of \cite{HZ}
breaks down when $p=p_c$ justifying this way our new paper. In \cite{HZ}, we noticed that for the unperturbed equation (\ref{mu=0})  with $p\le p_c$, Merle and Zaag \cite {MZ1}, \cite{MZ2} and  \cite{MZ3} crucially rely on the existence of a Lyapunov functional  in similarity variables established by Antonini and Merle \cite{AM}. Following this idea, we introduce  in \cite{HZ} (for $p\le p_c$)  the similarity variables   defined,
for all $x_0\in \er^N$,    $0< T_0\le T(x_0)$ by
\begin{equation}\label{scaling}
y=\frac{x-x_0}{T_0-t},\qquad s=-\log (T_0-t),\qquad u(x,t)=\frac{1}{(T_0-t)
^{\frac{2}{p-1}}}w_{x_0,T_0}(y,s).
\end{equation}
From (\ref{1}), the  function $w_{x_0,T_0}$  (we write $w$ for simplicity)
satisfies the following equation for all $y\in B\equiv B(0,1)$ and $s\ge -\log T_0$:
\begin{eqnarray}\label{B}
\partial_{ss}w&\! \!\!=\! \!\!&\frac{1}{\rho_{\alpha}}\div(\rho_{\alpha} \grad w-\rho_{\alpha}(y.\grad w)y)
-\frac{2p+2}{(p-1)^2}w+|w|^{p-1}w
-\frac{p+3}{p-1}\partial_sw-2y.\grad \partial_sw\nonumber\\
&&+e^{-\frac{2ps}{p-1}}f\Big(e^{\frac{2s}{p-1}}w\Big)+
e^{-\frac{2ps}{p-1}}g\Big(e^{\frac{(p+1)s}{p-1}}(\partial_sw+y.\grad w+\frac{2}{p-1}w)\Big),\qquad\quad
\end{eqnarray}
where $\rho_{\alpha}=(1-|y|^2)^{\alpha}$, with $\alpha=\alpha(N,p)=\frac{2}{p-1}-\frac{N-1}2$.
In the new set of variables $(y,s),$
the behavior of $u$ as $t \rightarrow T_0$
is equivalent to the
behavior of $w$ as $s \rightarrow +\infty$.

\bigskip

\no Following  Antonini and Merle \cite{AM} and  Merle and Zaag  \cite{MZ2}, we multiply equation (\ref{B}) by $\rho_{\alpha}w_s$ and integrate over the unit ball $B$.

\bigskip

\no When $(f,g)\equiv(0,0)$, we readily see from this calculation, as in  \cite{AM}  and   \cite{MZ2}, that $E_{\alpha}$ (defined below in (\ref{energyf})) is a Lyapunov functional  in the sense that
\begin{equation}\label{I1}
\frac{d}{ds}E_{\alpha}(w)=
-2\alpha \int_{B}(\partial_sw)^2\frac{\rho_{\alpha}}{1-|y|^2}{\mathrm{d}}y,\qquad  ({\textrm{if}} \ p<p_c).\qquad
\end{equation}
\begin{equation}\label{I2}
\frac{d}{ds}E_0(w)=
- \int_{\partial B}(\partial_sw)^2{\mathrm{d}}\sigma,\qquad  ({\textrm{if}} \ p=p_c).\qquad\qquad\qquad
\end{equation}
\bigskip

\no
When $(f,g)\not\equiv(0,0)$, (\ref{I1}) and (\ref{I2}) are perturbed by exponentially small terms with no sign. Our idea in  \cite{HZ} was to control these perturbation terms by the terms appearing  in the definition of $E_{\alpha}(w)$  or its dissipation given in (\ref{I1}) and (\ref{I2}) where
$(f,g)=(0,0)$. Since these perturbations terms are supported in the whole unit ball $B$  for $p\le p_c$, our study works when $p<p_c$, because the dissipation in
(\ref{I1}) is also supported in $B$. When $p=p_c$, the dissipation in (\ref{I2})
degenerates to the boundary and the method of \cite{HZ}
breaks down.
That obstruction fully justifies our new paper, where we invent a new idea to get a Lyapunov functional for equation when $p=p_c$ (note
the that $\alpha(N,p_c)=0$).
 Let us explain in the following our new idea. In fact, our strategy relies are two steps:

\bigskip
{\bf{Part 1: A rough estimate}}. As we said above, if we study $E_0(w)$, the functional of the unperturbed case $(f,g)\equiv (0,0)$, then we
get exponentially small terms that we can't control (unlike the subcritical case $p<p_c$). In other words, the study of $E_0(w)$ cannot be extended from the  unperturbed case to the general case. Fortunately, we noticed that the study of   $E_{\eta}(w)$ where $\eta>0$ can be extended from
the case $(f,g)\equiv (0,0)$ to the case $(f,g)\not\equiv (0,0)$.
It happens that in the former case,  $E_{\eta}(w)$ has a very bad bound, in the sense that $$E_{\eta}(w(s))\le Ce^{\beta \eta s},$$ for some $\beta=\beta (p)$ and
\begin{equation}\label{rough}
\|w(s)\|_{H^1(B)}+\|\partial_s w(s)\|_{L^2(B)}\le Ce^{\beta \eta  s}.
\end{equation}
 This behavior is conserved when $(f,g)\not\equiv (0,0)$.

\bigskip
{\bf{Part 2: The sharp  estimate}}. Now, we go back to $E_0(w)$, which is
the ''good'' functional for $(f,g)\equiv (0,0)$, in the sense that it is bounded. When $(f,g)\not\equiv (0,0)$, if we study
$E_0(w)$, the exponentially small terms are no longer a problem for us, thanks the rough estimate of (\ref{rough}), provided that we fix $\eta$ small enough.

\bigskip

\no
The equation (\ref{B}) will be studied in the space $\cal H$
$${\cal H}=\Big \{(w_1,w_2), |
\displaystyle\int_{B}\Big ( w_2^2
+|\grad w_1|^2(1-|y|^2)+w_1^2\Big ) {\mathrm{d}}y<+\infty \Big \}.$$
In the whole paper,   we  denote
$\ds{F(u)=\int_0^uf(v)dv}$
and we assume that (\ref{pc}) holds.

\bigskip

\no
In the case $(f,g)\equiv(0,0)$,   Merle  and Zaag \cite{MZ2} proved that
\begin{eqnarray}
\label{f}
E_0(w)=\displaystyle\int_{B}\Big ( \frac{1}{2}(\partial_sw)^2
+\frac{1}{2}|\grad w|^2-\frac{1}{2}(y.\grad w)^2+\frac{p+1}{(p-1)^2}w^2
-\frac{1}{p+1}|w|^{p+1}\Big ) {\mathrm{d}}y,\quad
\end{eqnarray}
is a Lyapunov functional for equation (\ref{B}).
 When $(f,g) \not\equiv (0,0)$, we introduce
\begin{eqnarray}
\label{f1}
 H(w)&=&E(w)+\sigma e^{-\gamma s},
\end{eqnarray}
where $\sigma$  is a sufficiently large constant that will be determined later,
\begin{eqnarray}
\label{f2}
 E(w)&=&E_0(w)+I_0(w)\quad  {\textrm {and}}\qquad
  I_0(w)=- e^{-\frac{2(p+1)s}{p-1}}\displaystyle\int_{B}F(e^{\frac{2s}{p-1}}w) {\mathrm{d}}y,\nonumber\\
   &&\qquad {\textrm {with}} \qquad \gamma=\min(\frac12,\frac{p-q}{p-1}
 )>0.
\end{eqnarray}

\bigskip
\no
Here we announce our main result.

\begin{theo}(Existence  of a Lyapunov functional for equation ({\ref{B}}))\\
\label{lyap}
 Consider
 $u $   a solution of ({\ref{1}}) with blow-up graph $\Gamma:\{x\mapsto T(x)\}$ and  $x_0$ is a
non characteristic point. Then
  there exists  $t_0(x_0)\in [0, T(x_0))$
such that,  for all $T_0\in
 (t_0(x_0), T_0(x_0)]$, for all  $s\ge -\log (T_0-t_0(x_0))$,     we have
\begin{eqnarray}\label{L1}
H(w(s+10))-H(w(s))
&\le&  - \int_{s}^{s+10}\int_{\partial B}(\partial_sw)^2(\sigma,\tau)d\sigma d\tau, \qquad \qquad
\end{eqnarray}
where $w=w_{x_0,T_0}$ is defined in (\ref{scaling}).
\end{theo}
\no {\bf{Remark 1.1.}}

\no $\bullet$ Since we crucially need 
a covering technique in our argument (see Appendix A), in fact, we need
a uniform version for  $x$ near $x_0$ (see theorem 1.1' page
\pageref{lyapbis} below).\\

\no $\bullet$ One may wander why  we take only sublinear perturbations in $\partial_tu$ (see hypothesis ($H_g$)).
It happens that  any superlinear  terms in $(\partial_tu $ generates in similarity variables $L^r$
norms of $\partial_sw$ and $\grad w$, where $r>2$, hence,
non controllable by the terms in the Lyapunov functional $E_0(w)$ (\ref{f}) of the non perturbed equation (\ref{mu=0}).

\bigskip

\no
The existence of this Lyapunov functional (and a blow-up criterion for equation (\ref{B})
based in $H$)  are a crucial step in the derivation of the blow-up rate for equation (\ref{1}).
Indeed, with the  functional $H$ and some more work, we are able to adapt the analysis
performed in \cite{MZ3} for equation (\ref{mu=0}) and get the following result:
\begin{theo}\label{t}(Blow-up rate  for equation ({\ref{1}}))\\
There exist
 $\varepsilon_0>0$   such that if $u $   is a solution of ({\ref{1}})
with blow-up graph $\Gamma:\{x\mapsto T(x)\}$ and  $x_0$ is a non characteristic point, then
there exist  $\widehat{S}_2>0$ such that

i)
 For all
 $s\ge \widehat{s}_2(x_0)=\max(\widehat{S}_2,-\log \frac{T(x_0)}4)$,
\begin{equation*}
0<\varepsilon_0\le \|w_{x_0,T(x_0)}(s)\|_{H^{1}(B)}+
\|\partial_s w_{x_0,T(x_0)}(s)\|_{L^{2}(B)}
\le K,
\end{equation*}
where $w_{x_0,T(x_0)}$ is defined in (\ref{scaling}) and $B$  is the unit ball of $\er^N$.\\
ii)  For all
  $t\in [t_2(x_0),T(x_0))$, where  $t_2(x_0)=T(x_0)-e^{-\widehat{s}_2(x_0)}$, we have
\begin{eqnarray*}
&&0<\varepsilon_0\le (T(x_0)-t)^{\frac{2}{p-1}}\frac{\|u(t)\|_{L^2(B(x_0,{T(x_0)-t}))}}{ (T(x_0)-t)^{\frac{N}{2}}}\nonumber\\
&&+ (T(x_0)-t)^{\frac{2}{p-1}+1}\Big (\frac{\|\partial_tu(t)\|_{L^2(B(x_0,{T(x_0)-t}))}}{ (T(x_0)-t)^{\frac{N}{2}}}+
 \frac{\|\grad u(t)\|_{L^2(B(x_0,{T(x_0)-t}))}}{ (T(x_0)-t)^{\frac{N}{2}}}\Big )\le K,
\end{eqnarray*}
where  $K=K(\widehat{s}_2(x_0),\|(u(t_2(x_0)),\partial_tu(t_2(x_0)))\|_{
H^{1}\times L^{2}(B(x_0,\frac{e^{-\widehat{s}_2(x_0)}}{\delta_0(x_0)})
)})$ and $\delta_0(x_0)\in (0,1)$ is defined in (\ref{cone1}).
\end{theo}
\no {\bf{Remark 1.2.}} With this blow-up rate, one can ask whether the results proved by Merle and
Zaag for the non perturbed problem in \cite{MZ4} \cite{MZ5} \cite{MZarxiv}, hold
for equation (\ref{1}) (blow-up, profile, regularity of the blow-up graph, existence
of characteristic points, etc...).  We believe that it is the case, however, the proof
will be highly technical, with no interesting ideas (in particular, equation
(\ref{1}) is not conserved under the Lorentz transform, which is crucial in  \cite{MZ4}
 \cite{MZ5} \cite{MZarxiv}, and lots of minor term will
appear in the analysis). Once again, we believe that the key point in the analysis of blow-up for
equation (\ref{1}) is the derivation of a Lyapunov functional
in similarity variables, which is the object of our paper.

\medskip

\no
As in the particular case where $(f,g)\equiv (0,0)$, the proof of theorem \ref{t} relies on four ideas (the existence of a Lyapunov functional,  interpolation in Sobolev spaces, some
critical Gagliardo-Nirenberg estimates and
a covering technique adapted to the geometric shape of the blow-up surface). It happens that adapting the proof of
\cite{MZ3} given in the non perturbed case (\ref{mu=0}) is straightforward, except for a key
argument, where we bound the $L^{p+1}$ space-time norm of $w$. Therefore, we only present that argument, and refer to \cite{HZ},
\cite{MZ1} and \cite{MZ3}  for the rest of the proof.


\medskip

\noindent This paper is organized as follows: In section 2, we obtain a rough control on space-time of the solution $w$.
Using this result, we proved in section 3, that the ''natural'' functional is a Lyapunov functional for equation
 (\ref{B}). Thus, we get theorem \ref{lyap}. Finaly,
 applying  this last theorem and  method used in \cite{MZ3}, we easily
 prove theorem~\ref{t}.

\section{A rough estimate}
In this section, we prove a rough version of $i)$ of  theorem \ref{t}, where we obtain an exponentially growing bound
on time averages of the $H^1\times L^2(B)$ norm of $(w,\partial_sw)$.
Consider
 $u $  a solution of ({\ref{1}}) with blow-up graph $\Gamma:\{x\mapsto T(x)\}$ and  $x_0$ is a non characteristic point.
 More precisely, this is the aim of this section.
\begin{pro}
\label{tt}
For all $\eta \in (0,1)$,
  there exists  $t_0(x_0)\in [0, T(x_0))$
such that,  for all $T_0\in
 (t_0(x_0), T(x_0)]$, for all $s\ge -\log (T_0-t_0(x_0))$ and
$x\in \er^N$ where $|x-x_0|\le \frac{e^{-s}}{\delta_0(x_0)}$,
   we have
\begin{eqnarray}\label{cor0004}
\int_{s}^{s+10} \!\int_{B}\!(\partial_s w)^2(y,\tau){\mathrm{d}}y{\mathrm{d}}\tau
+\int_{s}^{s+10} \int_{B}|w(y,\tau)|^{p+1}{\mathrm{d}}y{\mathrm{d}}\tau\nonumber\\
 + \int_{s}^{s+10} \int_{B}|\grad w(y,\tau)|^2{\mathrm{d}}y{\mathrm{d}}\tau\le K_1e^{\frac{\eta (p+3)s}{2}},
\end{eqnarray}
where $w=w_{x,T^*(x)}$  is defined in (\ref{scaling}) with
\begin{equation}\label{explosionx}
T^*(x)=T_0-\delta_0(x_0)(x-x_0),
\end{equation}
$K_1=K_1(\eta,T(x_0)-t_0(x_0),\|(u(t_0(x_0)),\partial_tu(t_0(x_0)))\|_{
H^{1}\times L^{2}(B(x_0,\frac{T(x_0)-t_0(x_0)}{\delta_0(x_0)})
)})$
and
$\delta_0(x_0)\in (0,1)$ is defined in (\ref{cone1}).
\end{pro}

\subsection{Rough energy estimates}
Consider   $T_0\in (0, T(x_0)]$, for all 
$x\in \er^N$ is such that $|x-x_0|\le \frac{T_0}{\delta_0(x_0)}$, where $\delta(x_0)$ is defined in (\ref{cone1}), then we write $w$ instead of $w_{x,T^*(x)}$ defined in (\ref{scaling}) with $T^*(x)$ given in  (\ref{explosionx}).
Let  $\eta \in (0,1)$ and
 write the equation
(\ref{B}) satisfied by $w$ in the form
\begin{eqnarray}\label{C}
\partial_{ss}w&=&\frac{1}{\rho_{\eta}}\div(\rho_{\eta} \grad w-\rho_{\eta}(y.\grad w)y)+2\eta
 y.\grad w
-\frac{2p+2}{(p-1)^2}w+|w|^{p-1}w\nonumber\\
&&-\frac{p+3}{p-1}\partial_s w-2y.\grad \partial_sw
+e^{-\frac{2ps}{p-1}}f\Big(e^{\frac{2s}{p-1}}w\Big)\\ &&+
e^{-\frac{2ps}{p-1}}g\Big(e^{\frac{(p+1)s}{p-1}}(\partial_s w+y.\grad w+\frac{2}{p-1}w)\Big),\
  \forall y\in B\ {\textrm{ and}} \ s\ge -\log T^*(x),\qquad\quad\nonumber
\end{eqnarray}
where $\rho_{\eta}=(1-|y|^2)^{\eta}$.
We denote by $C$ a constant which depends on $\m$.

\medskip

\noindent
To control the norm of $(w(s),\partial_s w(s))\in {\cal H}$, we first introduce the following functionals
\begin{eqnarray}
\label{energyf}
E_{\eta}(w)&=&\displaystyle\int_{B}\Big ( \frac{1}{2}(\partial_s w)^2
+\frac{1}{2}|\grad w|^2-\frac{1}{2}(y.\grad w)^2+\frac{p+1}{(p-1)^2}w^2
-\frac{1}{p+1}|w|^{p+1}\Big )\rho_{\eta} {\mathrm{d}}y,\nonumber\qquad\\
  I_\eta(w)&=&- e^{-\frac{2(p+1)s}{p-1}}\displaystyle\int_{B}F(e^{\frac{2s}{p-1}}w)\rho_{\eta} {\mathrm{d}}y,\nonumber\\
 J_\eta(w)&=& -\eta\displaystyle\int_{B}w\partial_s w\rho_{\eta} {\mathrm{d}}y+ \frac{N\eta}2 \displaystyle\int_{B}w^2 \rho_{\eta}{\mathrm{d}}y,\\
 H_\eta(w)&=&\displaystyle {E_\eta(w)+I_\eta(w)+J_\eta(w),}\nonumber\quad\\
 G_\eta(w)&=&\displaystyle H_\eta(w)\ e^{-\frac{\eta(p+3)s}{2}}+\theta e^{-\frac{\eta(p+3)s}{2}},\nonumber\quad
\end{eqnarray}
where $\theta=\theta (\eta)$  is a sufficiently large constant that will be determined later.
In this subsection, we prove that
$ G_\eta(w)$ is decreasing in time, which will  give the rough (ie exponentially fast)
estimate for  $E_\eta(w)$ and $\|(w,\partial_s w)\|_{{H^1(B)}\times{L^2(B)}}$.

\no
Now we state two lemmas which are crucial for the proof. We begin with   bounding the time derivative of
$ E_\eta(w)$ in the following lemma.
\begin{lem}
\label{energylyap0}
For all  $\eta \in (0,1)$,   we have
 the following inequality, for all $s\ge \max (-\log T^*(x),0)$,
\begin{eqnarray}\label{energylem21}
\frac{d}{ds}(E_{\eta}(w)+I_{\m} (w))\!\le\!
-2\m\int_{B}(\partial_s w)^2\frac{|y|^2\rho_{\eta}}{1-|y|^2}{\mathrm{d}}y
+2\m \int_{B} \partial_s w (y.\grad w)\rho_{\eta}{\mathrm{d}}y
+\Sigma_ {0}(s),\qquad
\end{eqnarray}
where $\Sigma_0(s)$ satisfies
\begin{eqnarray}
\label{energysigm}
\Sigma_ {0}(s) &\le&\!\!\!\!Ce^{-2\gamma s}+Ce^{-2\gamma s}\int_{B}|\grad w|^2(1-|y|^2)\rho_{\eta}{\mathrm{d}}y+Ce^{-2\gamma s}\int_{B}w^2\rho_{\eta}{\mathrm{d}}y\nonumber\\
&&+Ce^{-2\gamma s}\int_{B}(\partial_s w)^2  \frac{\rho_{\eta}}{1-|y|^2}{\mathrm{d}}y
+Ce^{-2\gamma s}\int_{B}\!\!|w|^{p+1}\rho_{\eta} {\mathrm{d}}y,\qquad \qquad
\end{eqnarray}
with  $\gamma=\min(\frac12, \frac{p-q}{p-1} )>0.$
\end{lem}
{\it{Proof}}:  Multipling $(\ref{C})$ by $\rho_{\eta}\partial_s w $ and integrating over the ball $B$, we obtain for all $s\ge -\log T^*(x)$,
\begin{eqnarray}\label{energyE0}
\frac{d}{ds}(E_{\eta}(w)+I_{\m} (w))&=&
-2 \int_{B} (\partial_s w)(y.\grad \partial_s w)\rho_{\eta}{\mathrm{d}}y
-\frac{p+3}{p-1} \int_{B}(\partial_s w)^2\rho_{\eta}{\mathrm{d}}y\nonumber\\
&&+2\m \int_{B} (\partial_s w)(y.\grad w)\rho_{\eta}{\mathrm{d}}y\nonumber\\
&&+\underbrace{\frac{2(p+1)}{p-1}e^{-\frac{2(p+1)s}{p-1}}\int_{B}F\Big(e^{\frac{2s}{p-1}}w\Big)\rho_{\eta} {\mathrm{d}}y}_{\Sigma_0^1(s)} \nonumber\\
&&+\underbrace{\frac{2}{p-1}e^{-\frac{2ps}{p-1}}\int_{B}f\Big(e^{\frac{2s}{p-1}}w\Big)w\rho_{\eta} {\mathrm{d}}y }_{\Sigma_0^2(s)}\\
&&+\underbrace{ e^{-\frac{2ps}{p-1}}\int_{B}g\Big(e^{\frac{(p+1)s}{p-1}}(\partial_s w+y.\grad w+\frac{2}{p-1}w)\Big)\partial_s w\rho_{\eta} {\mathrm{d}}y}_{\Sigma_0^3(s)}.\nonumber
\end{eqnarray}
Then by integrating by parts, we have
\begin{eqnarray}\label{energym}
\frac{d}{ds}(E_{\eta}(w)+I_{\m} (w))&=&
-2 \m \int_{B} (\partial_s w)^2 \frac{|y|^2\rho_{\eta}}{1-|y|^2}{\mathrm{d}}y
+2\m \int_{B} (\partial_s w)(y.\grad w)\rho_{\eta}{\mathrm{d}}y\nonumber\\
&&+\Sigma_0^1(s) +\Sigma_0^2(s)+\Sigma_0^3(s).\quad
\end{eqnarray}
Using the fact that
$ |{F(x)}|+|x{f(x)}|\le  C( 1+|x|^{p+1}),$ we obtain that for all $s \ge \max (-\log T^*(x),0)$,
\begin{eqnarray}\label{energyIm}
|\Sigma_0^1(s)|+|\Sigma_0^2(s)|&\le& Ce^{-2\gamma s}+ Ce^{-2\gamma s}\int_{B}|w|^{p+1}\rho_{\eta} {\mathrm{d}}y,
\end{eqnarray}
on the one hand. On the other hand,  since
$ |g(x)|\le M( 1+|x|)$, we write that
for all $s \ge \max (-\log T^*(x),0)$,
\begin{eqnarray*}\label{energyIm1}
|\Sigma_0^3(s)|&\le& Ce^{-2\gamma s}
\int_{B}(\partial_s w)^2\rho_{\eta} {\mathrm{d}}y +
Ce^{-2\gamma s}\int_{B}|y.\grad w ||\partial_s w|\rho_{\eta} {\mathrm{d}}y\nonumber\\
&&+ C e^{-2\gamma s}\int_{B}|w\partial_s w|\rho_{\eta} {\mathrm{d}}y+ C e^{-2\gamma s}\int_{B}|\partial_s w|\rho_{\eta} {\mathrm{d}}y.\qquad
\end{eqnarray*}
By exploiting the inequality $ab\le \frac{a^2}{2}+\frac{b^2}{2},$  we conclude that
for all \\$s \ge \max (-\log T^*(x),0),$
\begin{eqnarray}\label{energyIm2}
|\Sigma_0^3(s)|&\le& Ce^{-2\gamma s}\int_{B}(\partial_s w)^2\frac{\rho_{\eta}}{1-|y|^2}{\mathrm{d}}y+C  e^{-2\gamma s}\int_{B}w^2\rho_{\eta}{\mathrm{d}}y\nonumber\\
&&+C e^{-2\gamma s}\int_{B}|\grad w|^2(1-|y|^2)\rho_{\eta}{\mathrm{d}}y+ C e^{-2\gamma s}.
\end{eqnarray}
Then, by using  (\ref{energym}), (\ref{energyIm}) and (\ref{energyIm2}), we have
for all $s \ge \max (-\log T^*(x),0)$, the estimates (\ref{energylem21}) and  (\ref{energysigm}) hold.
 This concludes the proof of lemma
 \ref{energylyap0}.

\Box

\vspace{0.3cm}

\no We are now going to prove the following estimate for the functional $J_{\m}$:
\begin{lem}
\label{energylyap1}
For all $\eta \in (0,1)$,  $J_{\m}$ satisfies
 the following inequality, for all $s\ge \max (-\log T^*(x),0)$
\begin{eqnarray}\label{energytheta1}
\frac{d}{ds}J_{\m}(w)&\le& \frac{32\m}{p+15}\!\int_{B}\!(\partial_s w)^2\frac{\rho_{\eta}}{1-|y|^2}{\mathrm{d}}y-2\m
\int_{B}\partial_s w(y.\grad w) \rho_{\eta}{\mathrm{d}}y+\frac{\m (p+3)}{2}H_{\m}\nonumber\\ &&
-\frac{\m (p+15)}{8} \int_{B}(\partial_s w)^2\rho_{\eta}{\mathrm{d}}y
-  \frac{\m (p-1)}{2(p+1)}\int_{B}|w|^{p+1}\rho_{\eta}{\mathrm{d}}y\\ &&
- \frac{ \m( p-1)}{8}\int_{B}|\grad w|^2(1-|y|^2)\rho_{\eta}{\mathrm{d}}y+\Sigma_ {1}(s),\nonumber
\end{eqnarray}
where   $\Sigma_1(s)$ satisfies
\begin{eqnarray}
\label{energysigm1}
 \Sigma_ 1(s) &\le&
C  e^{-2\gamma s}\!\int_{B}\!(\partial_s w)^2\frac{\rho_{\eta}}{1-|y|^2}{\mathrm{d}}y
+C   e^{-2\gamma s}\int_{B}|w|^{p+1}\rho_{\eta}{\mathrm{d}}y
\\ &&+C  e^{-2\gamma s}\int_{B}|\grad w|^2(1-|y|^2)\rho_{\eta}{\mathrm{d}}y+C \int_{B}w^2 \rho_{\eta}{\mathrm{d}}y+C  e^{-2\gamma s},\nonumber
\end{eqnarray}
with $\gamma=\min(\frac12,\frac{p-q}{p-1})$.
\end{lem}
{\it Proof:}
Note that  $J_{\m}$ is a differentiable  function for  all
 $s\ge   -\log T^*(x)$ and that
\begin{eqnarray*}
\frac{d}{ds}J_{\m}(w)&=&- \m\int_{B}(\partial_sw)^2\rho_{\eta}{\mathrm{d}}y- \m\int_{B}w\partial_{ss}w\rho_{\eta}{\mathrm{d}}y
+N\m \int_{B}w\partial_s w \rho_{\eta}{\mathrm{d}}y.
\end{eqnarray*}
By using  equation (\ref{C}) and integrating by parts, we have
\begin{eqnarray}\label{energyt1}
\frac{d}{ds}J_{\m}(w)&=&
-\eta \int_{B}(\partial_s w)^2\rho_{\eta}{\mathrm{d}}y
+ \m\int_{B}(|\grad w|^2-(y.\grad w)^2)\rho_{\eta}{\mathrm{d}}y- \m\int_{B}|w|^{p+1}\rho_{\eta}{\mathrm{d}}y\nonumber\\
&&-2\m \int_{B}\partial_s w(y.\grad w) \rho_{\eta}{\mathrm{d}}y +\m (\frac{2p+2}{(p-1)^2}+\m N)\int_{B}w^2\rho_{\eta}{\mathrm{d}}y\nonumber\\
&&+\underbrace{4\m^2 \int_{B}w\partial_s w\frac{|y|^2\rho_{\eta}}{1-|y|^2} {\mathrm{d}}y}_{\Sigma_1^1(s)}
\underbrace{-2\m^3\int_{B}w^2\frac{|y|^2\rho_{\eta}}{1-|y|^2} {\mathrm{d}}y}_{\Sigma_1^2(s)}\\
&&\underbrace{-\m e^{-\frac{2ps}{p-1}}\int_{B}wf\Big(e^{\frac{2s}{p-1}}w\Big){\rho_{\eta}}{\mathrm{d}}y}_{\Sigma_1^3(s)}\nonumber\\
&&\underbrace{ -\m e^{-\frac{2ps}{p-1}}\int_{B}wg\Big(e^{\frac{(p+1)s}{p-1}}(\partial_s w+y.\grad w+\frac{2}{p-1}w)\Big){\rho_{\eta}}{\mathrm{d}}y}_{\Sigma_1^4(s)}.
\nonumber
\end{eqnarray}
By combining (\ref{energyf}) and (\ref{energyt1}),   we write
\begin{eqnarray}\label{energytheta}
\frac{d}{ds}J_{\m}(w)&=&-2\m \int_{B}\partial_s w(y.\grad w) \rho_{\eta}{\mathrm{d}}y+\frac{\m (p+3)}{2}H_{\m}
-\frac{\m (p+7)}{4} \int_{B}(\partial_s w)^2\rho_{\eta}{\mathrm{d}}y\nonumber\\ &&- \frac{\m(p-1)}{2(p+1)}\int_{B}|w|^{p+1}\rho_{\eta}{\mathrm{d}}y
-\frac{ \m (p-1)}{4}\int_{B}(|\grad w|^2-(y.\grad w)^2)\rho_{\eta}{\mathrm{d}}y\\
&&-\m \Big(\frac{ p+1}{2(p-1)}+\frac{  \m N( p-1)}{4} \Big)  \int_{B}w^2 \rho_{\eta}{\mathrm{d}}y
+ \underbrace{\frac{\m^2(p+3) }{2}\displaystyle\int_{B}w\partial_s w\rho_{\eta} {\mathrm{d}}y}_{\Sigma_1^5(s)}\nonumber\\
&&+\underbrace{\frac{\m (p+3)}2  e^{-\frac{2(p+1)s}{p-1}}\displaystyle\int_{B}F(e^{\frac{2}{p-1}s}w)\rho_{\eta} {\mathrm{d}}y}_{\Sigma_1^6(s)}+\Sigma_1^1(s)+\Sigma_1^2(s)+\Sigma_1^3(s)+\Sigma_1^4(s).\nonumber
\end{eqnarray}
We now study each of the last six terms.

\no
By the Cauchy-Schwartz inequality we write, for all $\mu\in (0,1)$
\begin{eqnarray}\label{energyJ2}
\Sigma_1^1(s)
&\le&2\m (1-\mu)\int_{B}(\partial_s w)^2\frac{\rho_{\eta}}{1-|y|^2}{\mathrm{d}}y+\frac{2 \m^3}{1-\mu}\int_{B}w^2\frac{|y|^2\rho_{\eta}}{1-|y|^2}
{\mathrm{d}}y.
\qquad\nonumber
\end{eqnarray}
We use the expression of $\Sigma_1^2(s)$ to obtain, for all $ \mu\in (0,1)$
\begin{eqnarray}\label{energyJ2J3}
\Sigma_1^1(s)+\Sigma_1^2(s)
&\le&2\m (1-\mu)\int_{B}(\partial_s w)^2\frac{\rho_{\eta}}{1-|y|^2}{\mathrm{d}}y+\frac{2 \m^3\mu }{1-\mu}\int_{B}w^2\frac{|y|^2\rho_{\eta}}{1-|y|^2}
{\mathrm{d}}y.\qquad
\end{eqnarray}\no Since we have the following Hardy type inequality for any $w\in
H^{1}_{loc,u}(\er^N)$ (see appendix $B$  in \cite{MZ1} for details):
\begin{eqnarray}\label{energyhardyJJ}
\int_{B}w^2\frac{|y|^2\rho_{\eta}}{1-|y|^2}{\mathrm{d}}y
&\le&\frac{1}{\m^2}\int_{B}|\grad w|^2(1-|y|^2)\rho_{\eta}{\mathrm{d}}y+\frac{N}{\m} \int_{B}w^2\rho_{\eta}
{\mathrm{d}}y,
\end{eqnarray}
from (\ref{energyJ2J3}) and (\ref{energyhardyJJ}) and if  we choose
 $\mu =\frac{p-1 }{p+15}$,
 we conclude that
\begin{eqnarray}\label{energyJ22}
\Sigma_1^1(s)+\Sigma_1^2(s)&\le&\frac{32\m}{p+15} \!\int_{B}\!(\partial_s w)^2\frac{\rho_{\eta}}{1-|y|^2}{\mathrm{d}}y+
\frac{\m^2N(p-1)}{8} \int_{B}\!w^2\rho_{\eta}
{\mathrm{d}}y\nonumber\\
&&+\frac{\m(p-1)}{8} \int_{B}\!|\grad w|^2(1-|y|^2)\rho_{\eta}
{\mathrm{d}}y.\qquad
\end{eqnarray}
\no
By exploiting the fact that
$ |{F(x)}|+   |xf(x)|\le C(1+|x|^{p+1}),$
we obtain for all $s\ge \max(-\log T^*(x),0)$
\begin{eqnarray}\label{energyJ4}
\Sigma_1^3(s)+\Sigma_1^6(s)&\le&C  e^{-2\gamma s}+
 C e^{-2\gamma s}\int_{B}|w|^{p+1}\rho_{\eta}{\mathrm{d}}y.
\end{eqnarray}
\no
In a similar way,  by  using the fact that
$ |g(x)|\le M(1+|x|)$, we write for all ${s\ge \max(-\log T^*(x),0)}$
\begin{eqnarray*}\label{energyJ5}
\Sigma_1^4(s)&\le& C e^{-2\gamma s}
\int_{B}w^2\rho_{\eta} {\mathrm{d}}y +
C e^{-2\gamma s}\int_{B}|y.\grad w ||w|\rho_{\eta} {\mathrm{d}}y\nonumber\\
&&+ C e^{-2\gamma s}\int_{B}(\partial_s w)^2\rho_{\eta} {\mathrm{d}}y+ C  e^{-2\gamma s}.\qquad
\end{eqnarray*}
By using (\ref{energyhardyJJ}), we get
\begin{eqnarray}\label{energyJ55}
\Sigma_1^4(s)&\le& C e^{-2\gamma s}
\int_{B}(\partial_s w)^2\rho_{\eta} {\mathrm{d}}y +C e^{-2\gamma s}\int_{B}|\grad w|^2(1-|y|^2)\rho_{\eta}
{\mathrm{d}}y\nonumber\\
&& +C e^{-2\gamma s}\int_{B}w^2\rho_{\eta} {\mathrm{d}}y+ C e^{-2\gamma s}.\qquad
\end{eqnarray}
To estimate $\Sigma_1^5(s)$, we use  the Cauchy-Schwartz inequality and we get
\begin{eqnarray}\label{energyJ1}
\Sigma_1^5(s)& \le&
 \frac{(p-1)\m}{8} \int_{B}(\partial_s w)^2\rho_{\eta}{\mathrm{d}}y+
C  \int_{B}w^2\rho_{\eta}
{\mathrm{d}}y.\qquad\qquad
\end{eqnarray}
Since $|y.\grad w|\le |y||\grad w|$, it follows that
\begin{eqnarray}\label{energylast}
 \int_{B}|\grad w|^2(1-|y|^2)\rho_{\eta}{\mathrm{d}}y&\le &\int_{B}(|\grad w|^2-(y.\grad w)^2) \rho_{\eta}{\mathrm{d}}y.
\end{eqnarray}
Finally,
by using (\ref{energytheta}),  (\ref{energyJ22}),  (\ref{energyJ4}), (\ref{energyJ5}),
 (\ref{energyJ1}) and (\ref{energylast}),   we have easily  the estimate (\ref{energytheta1}) and (\ref{energysigm1}).
This concludes the proof of  lemma \ref{energylyap1}.
\Box

\bigskip
\no From lemmas \ref{energylyap0} and \ref{energylyap1}, we are in a position to prove the following   proposition
\begin{pro}(Existence  of  a decreasing  functional for equation ({\ref{B}}))\\
\label{lyapm}
For all $\eta \in (0,1)$, there exists  $S_0>0$
such that,  $G_{\m}$ defined in (\ref{energyf}) satisfies the following inequality, for all
$s_2>s_1\ge \max(-\log T^*(x),S_0)$,
\begin{eqnarray}\label{LE}
G_{\m}(w(s_2))-G_{\m}(w(s_1))
&\le&
- \frac{\m (p-1)}{p+15}\int_{s_1}^{s_2} e^{-\frac{\m (p+3)s}{2}}\!\int_{B}\!(\partial_s w)^2\frac{\rho_{\eta}}{1-|y|^2}{\mathrm{d}}y{\mathrm{d}}s\nonumber\\
&&-  \frac{\m (p-1)}{8(p+1)}\int_{s_1}^{s_2} e^{-\frac{\m (p+3)s}{2}}\int_{B}|w|^{p+1}\rho_{\eta}{\mathrm{d}}y{\mathrm{d}}s\\
&&- \frac{ \m( p-1)}{16} \int_{s_1}^{s_2} e^{-\frac{\m (p+3)s}{2}}\int_{B}|\grad w|^2(1-|y|^2)\rho_{\eta}{\mathrm{d}}y{\mathrm{d}}s.\nonumber
\end{eqnarray}
\end{pro}

\vspace{0.3cm}
{\it Proof:}
From lemmas \ref{energylyap0} and \ref{energylyap1}, we obtain for all $s\ge \max(-\log T^*(x), 0)$,
\begin{eqnarray}\label{m}
\frac{d}{ds}H_{\m}(w)&\le&\frac{\m (p+3)}{2}H_{\m}(w)
-\Big  (\frac{2\m(p-1)}{p+15} -C  e^{-2\gamma s}\Big )\!\int_{B}\!(\partial_s w)^2\frac{\rho_{\eta}}{1-|y|^2}{\mathrm{d}}y\nonumber\\
&&-  \Big (\frac{\m (p-1)}{2(p+1)}-C  e^{-2\gamma s}\Big)\int_{B}|w|^{p+1}\rho_{\eta}{\mathrm{d}}y\nonumber\\
&&- \Big (\frac{ \m( p-1)}{8}-C  e^{-2\gamma s}\Big )\int_{B}|\grad w|^2(1-|y|^2)\rho_{\eta}{\mathrm{d}}y\\
&&+C \int_{B}w^2 \rho_{\eta}{\mathrm{d}}y+C  e^{-2\gamma s}.\nonumber
\end{eqnarray}
We now choose $S_0$ large enough ($S_0\ge 0$), so that for all $s\ge S_0$, we have
\begin{eqnarray*}
 \frac{\m (p-1)}{p+15}-C e^{-2\gamma s} \ge 0, \qquad
 \frac{\m(p-1)}{4(p+1)}
-C e^{-2\gamma s} \ge 0,\quad \frac{\m(p-1)}{16}-C e^{-2\gamma s}\ge 0.
\end{eqnarray*}
Then, we deduce that, for all $s\ge \max(-\log T^*(x),S_0)$, we have
\begin{eqnarray}\label{m1}
\frac{d}{ds}H_{\m}(w)&\le&\frac{\m (p+3)}{2}H_{\m}(w)
-\frac{\m (p-1)}{p+15}\!\int_{B}\!(\partial_s w)^2\frac{\rho_{\eta}}{1-|y|^2}{\mathrm{d}}y-  \frac{\m (p-1)}{4(p+1)}\int_{B}|w|^{p+1}\rho_{\eta}{\mathrm{d}}y\nonumber\\
&&- \frac{ \m( p-1)}{16}\int_{B}|\grad w|^2(1-|y|^2)\rho_{\eta}{\mathrm{d}}y
+C \int_{B}w^2 \rho_{\eta}{\mathrm{d}}y+C  e^{-2\gamma s}.
\end{eqnarray}
By combining (\ref{m1}) and   the  following Jensen's inequality
\begin{eqnarray}\label{jj}
C\int_{B}w^2 \rho_{\eta}{\mathrm{d}}y\le \frac{\m (p-1)}{8(p+1)}
 \int_{B}|w|^{p+1}\rho_{\eta}{\mathrm{d}}y+C,
\end{eqnarray}
we obtain, for all $s\ge \max(-\log T^*(x),S_0)$,
\begin{eqnarray}\label{m2}
\frac{d}{ds}H_{\m}(w)&\le&\frac{\m (p+3)}{2}H_{\m}(w)
-\frac{\m (p-1)}{p+15}\!\int_{B}\!(\partial_s w)^2\frac{\rho_{\eta}}{1-|y|^2}{\mathrm{d}}y-  \frac{\m (p-1)}{8(p+1)}\int_{B}|w|^{p+1}\rho_{\eta}{\mathrm{d}}y\nonumber\\
&&- \frac{ \m( p-1)}{16}\int_{B}|\grad w|^2(1-|y|^2)\rho_{\eta}{\mathrm{d}}y
+C.
\end{eqnarray}
Finally, by (\ref{m2}), we prove easily that the
function $G_{\m}$ satisfies, for all \\ ${s\ge \max(-\log T^*(x),S_0),}$
\begin{eqnarray}\label{m3}
\frac{d}{ds}G_{\m}(w)\!\!&\le&\!\!\!
- \frac{\m (p-1)}{p+15}\ e^{-\frac{\m (p+3)s}{2}}\!\int_{B}\!\!(\partial_s w)^2\frac{\rho_{\eta}}{1-|y|^2}{\mathrm{d}}y-
\frac{\m (p-1)}{8(p+1)}\ e^{-\frac{\m (p+3)s}{2}}\!\int_{B}\!|w|^{p+1}\rho_{\eta}{\mathrm{d}}y\nonumber\\
&&\!\!- \frac{ \m( p-1)}{16}\ e^{-\frac{\m (p+3)s}{2}}\int_{B}\!|\grad w|^2(1-|y|^2)\rho_{\eta}{\mathrm{d}}y
+(C-\frac{\m \theta (p+3)}{2}) e^{-\frac{\m (p+3)s}{2}}.\nonumber
\end{eqnarray}
We now choose $\theta=\theta (\m)$ large enough,  so we have $C-\frac{\m \theta (p+3)}{2} \le 0$ and then
\begin{eqnarray}\label{m4}
\frac{d}{ds}G_{\m}(w)\!\!\!&\le&\!\!\!
-\frac{\m (p-1)}{p+15} \ e^{-\frac{\m (p+3)s}{2}}\!\int_{B}\!(\partial_s w)^2\frac{\rho_{\eta}}{1-|y|^2}{\mathrm{d}}y-  \frac{\m (p-1)}{8(p+1)}e^{-\frac{\m (p+3)s}{2}}\int_{B}\!|w|^{p+1}\rho_{\eta}{\mathrm{d}}y\nonumber\\
&&\!\!- \frac{ \m( p-1)}{16}e^{-\frac{\m (p+3)s}{2}}\int_{B}\!|\grad w|^2(1-|y|^2)\rho_{\eta}{\mathrm{d}}y.
\end{eqnarray}
Now (\ref{LE}) is a direct consequence of inequality
(\ref{m4}).
This concludes the proof of proposition \ref{lyapm}.
\Box


\subsection{Proof of proposition \ref{tt}}

\no We now claim the following lemma:
\begin{lem}\label{eL02}
For all $\eta \in (0,1)$, there exists $S_1\ge S_0$ such that, if
  $G_{\m}(w(s_1))<0$ for some $s_1\ge \max(-\log T^*(x),S_1)$, then $w$
blows up in some finite time $S>s_1$.
\end{lem}
{\it Proof:} The argument is the same as in the corresponding part
in \cite{HZ}. Let us remark   that our proof strongly relies on the fact that  $p<1+\frac{4}{N-2}$ which is implied by the fact that $p=p_c=1+\frac4{N-1}$.
\Box

\bigskip
\no We define the following time:
\begin{equation}\label{new1}
 t_0(x_0)=\max(T(x_0)-e^{-S_1},0).
\end{equation}

\no Since $\m\in (0,1)$, by combining  proposition \ref{lyapm} and  lemma \ref{eL02}, we get the following bounds:
\begin{cor}\label{cor10} (Estimates on $H_{\m}$)
For all $\eta \in (0,1)$,
  there exists  $t_0(x_0)\in [0, T(x_0))$
such that,  for all
 $T_0\in
 (t_0(x_0), T(x_0)]$, for all $s\ge -\log (T_0-t_0(x_0))$
 and $x\in \er^N$ where $|x-x_0|\le \frac{e^{-s}}{\delta_0(x_0)}$, we have

(i)
\begin{equation}\label{cor01}
-C\le H_{\m}(w(s))\le  \Big (\theta +  H_{\m}(w(\widetilde{s_0}))\Big )e^{\frac{\m(p+3)s}2},\nonumber
\end{equation}
\begin{equation}\label{corr01}
\int_{s}^{s+10} \!\int_{B}\!(\partial_s w)^2(y,\tau)\frac{\rho_{\eta}}{1-|y|^2}{\mathrm{d}}y{\mathrm{d}}\tau\le
C\big(\theta +  H_{\m}(w(\widetilde{s_0}))\big)e^{\frac{\m (p+3)s}{2}},
\end{equation}
\begin{eqnarray}\label{corrr01}
\int_{s}^{s+10}  \!\!\int_{B}\!|w(y,\tau)|^{p+1}\rho_{\eta}{\mathrm{d}}y{\mathrm{d}}\tau
 + \int_{s}^{s+10}\!\! \int_{B}\!|\grad w(y,\tau)|^2(1-|y|^2)\rho_{\eta}{\mathrm{d}}y{\mathrm{d}}\tau\nonumber\\
\le C\big(\theta +  H_{\m}(w(\widetilde{s_0}))\big) e^{\frac{\m (p+3)s}{2}}.\qquad\qquad
\end{eqnarray}
(ii)
\begin{equation}\label{cor02}
\int_{s}^{s+10} \!\!\int_{B}\!(\partial_s w)^2(y,\tau){\mathrm{d}}y{\mathrm{d}}\tau
\le C\big(\theta +  C H_{\m}(w(\widetilde{s_0}))\big)e^{\frac{\m (p+3)s}{2}},
\end{equation}
\begin{eqnarray}\label{cor002}
\int_{s}^{s+10} \int_{B_{\frac12}}|w(y,\tau)|^{p+1}{\mathrm{d}}y{\mathrm{d}}\tau
 + \int_{s}^{s+10} \int_{B_{\frac12}}|\grad w(y,\tau)|^2{\mathrm{d}}y{\mathrm{d}}\tau\nonumber\\
\le C\Big(\theta +  H_{\m}(w(\widetilde{s_0}))\Big)e^{\frac{\m (p+3)s}{2}},
\end{eqnarray}
where $w=w_{x,T^*(x)}$  is defined in (\ref{scaling}), $T^*(x)$ is defined in (\ref{explosionx}) and
$\widetilde{s_0}= -\log (T^*(x)-t_0(x_0))$.

\end{cor}
\no {\bf{Remark 2.2.}} By using the definition of  (\ref{scaling}) of
$w_{x,T^*(x)}=w$, we write easily
\begin{eqnarray}\label{cor0003}
 C\theta +  C H_{\m}(w(\widetilde{s_0}))
\le K_0,\nonumber
\end{eqnarray}
where
$K_0=K_0(\m,T(x_0)-t_0(x_0),\|(u(t_0(x_0)),\partial_tu(t_0(x_0)))\|_{
H^{1}\times L^{2}(B(x_0,\frac{T(x_0)-t_0(x_0)}{\delta_0(x_0)})
)})$
and
$\delta_0(x_0)\in (0,1)$ is defined in (\ref{cone1}).

\vspace{0.3cm}

\no {\it{Proof of proposition \ref{tt}:}}
In the following, we introduce a covering technique to derive  proposition \ref{tt} from  corollary \ref{cor10}.
  Note that the estimate on the space-time $L^2$ norm of $\partial_sw$ was already proved in (ii) of
 corollary \ref{cor10}. Thus, we focus on the  space-time $L^{p+1}$ norm of $w$ and  $L^2$ norm  $\grad w$.
 For that, we introduce a new covering technique to extend the
 estimate of any known space-time $ L^q$ norm of $w$, $\partial_s w$ or $\grad w$
 from $\B1$ to the whole unit ball. Note
 that, we follow the covering method of Merle and Zaag   in \cite{MZ3}. Therefore, we don't give all the details.
 Here, we strongly need the following
 local space-time generalization of the notion of characteristic point: a point $(x_0, T_0)\in \bar D_u$ is $\delta_0$-non characteristic with respect to $\t0$ where $\delta_0\in(0,1)$ if
\begin{equation}\label{nonchar2}
u \mbox{ is defined on }{\cal D}_{x_0, T_0,\t0, \delta_0}
\end{equation}
where
\begin{equation}\label{defP}
{\cal D}_{x_0,T_0, \t0, \delta_0}=\{(\xi,\tau)\neq (x_0, T_0)\;|\; \t0\le \tau\le T_0-\delta_0|\xi-x_0|\}.
\end{equation}
We also  define  ${\cal S}_{x_0,T_0, \t0, \delta_0}$ the slice of
${\cal D}_{x_0,T_0, \t0, \delta_0}$ betwen $\tau=t$ and $\tau=T_0-e^{-10} (T_0-t)$
 by
\begin{equation}\label{app2}
{\cal S}_{x_0,T_0, \t0, \delta_0}=\{(\xi,\tau)|\;
\t0\le \tau\le T_0-e^{-10} (T_0-t), |\xi-x_0|\le \frac{T_0-\tau}{\delta_0}\}.
\end{equation}
In fact, we find it easier to work in the $u(x,t)$ setting, in order to respect the geometry of the blow-up set.
We claim the following:

\begin{lem}[Covering technique]\label{covering00}Consider
 $\kappa\ge 0$, $q\ge 1$ and   $f\in L^q_{loc}\left( D_u\right)$.
 Then, for all $x_0\in \R^N$, $T_0\le T(x_0)$ and $t_1\le T_0$ such that
${\cal D}_{x_0, T_0,t_1, \delta_0}\subset D_u$ for some $\delta_0\in(0,1)$, we have:\\


\no (i) For any $x$ such that $|x-x_0|\le \frac{T_0-t_1}{\delta_0}$, $f$
is defined on  ${\cal S}_{x,T^*(x), t_1,1}$.
\[
(ii)\;\;\;\;\;\;\;\;\;\;\;\; \sup_{\{x\;|\;|x-x_0|\le
 \frac{T_0-t_1}{\delta_0}\}}\int_{
 {\cal S}_
 {x,T^*(x),t_1,1}} T^*(x)-\t0)^\kappa
  |f(\xi,\t0)|^q d\xi dt< +\infty.
\;\;\;\;\;\;\;\;\;\;\;\;\;\;\;\;\;\;\qquad
\]\qquad
\begin{eqnarray*}
(iii)\;\;\;\;\;\;\;\;\;\;&&\sup_{\{x\;|\;|x-x_0|\le
 \frac{T_0-t_1}{\delta_0}\}}
\int_{t_1}^{t_2(x)}
(T^*(x)-\t0)^\kappa
\int_{B(x, T^*(x)-t)} |f(\xi,\t0)|^q d\xi dt\;\;\\
&\le& C(\delta_0,\kappa)\sup_{\{x\;|\;|x-x_0|\le \frac{T_0-t_1}{\delta_0}\}}\int_{t_1}^{t_2(x)} (T^*(x)-t)^\kappa\int_{B(x, \frac{T^*(x)-\t0}2)} |f(\xi,\t0)|^q d\xi dt,
\end{eqnarray*}
where $t_2(x)=T^*(x)-e^{-10}(T^*(x)-t_1)$ and where $T^*(x)$ is defined in (\ref{explosionx}).

\end{lem}
\no {\bf Remarks 2.2.}

$\bullet$ The point $(x, T^*(x))$ is on the lateral boundary of ${\cal D}_{x_0, T_0, \t0, \delta_0}$.\\

$\bullet$ Note that the supremum is taken over the basis of ${\cal D}_{x_0, T_0,t_1, \delta_0}$.\\

\no {\it Proof of lemma \ref{covering00}: } See Appendix A.
 \Box\\

In the following, we give the proof only for the space-time $L^{p+1}$ norm of $w$, since the space-time $L^2$ norm of $\grad w$ follows exactly in the same way. Consider
  $x_0$ is a non characteristic point and $T_0\in (t_0(x_0), T(x_0)]$, where $t_0(x_0)$ is defined in (\ref{new1}). Consider  $s\ge -\log (T_0-t_0(x_0))$, consider then
$x\in \er^N$ such that, $|x-x_0|\le \frac{e^{-s}}{\delta_0(x_0)}$. If $w=w_{x,T_0}$, then we write from the
self-similar change of variables (\ref{scaling})
\begin{eqnarray}\label{b1}
\int_{s}^{s+10} \!\!\int_{B(0,\rho)}\!|w_{x,T^*\!(x)}(y,\tau)|^{p+1}{\mathrm{d}}y{\mathrm{d}}\tau=
\int_{t_1}^{t_2(x)}\!\!
\int_{B(x, \rho(T^*(x)-t))} |u(\xi,\t0)|^{p+1} d\xi dt,
\end{eqnarray}
 where $t_1=t_1(s,T_0)=T_0-e^{-s}$, $t_2(x)=t_2(x,s,T_0)=T^*(x)-e^{-10}(T^*(x)-t_1)$,  with $\rho=1$ or $\frac12$.
Note that
${\cal D}_{x_0, T_0,t_1, \delta_0(x_0)}\subset {\cal D}_{x_0, T(x_0),0, \delta_0(x_0)}\subset
{\cal D}_{u}$ by definition (\ref{cone1}). Therefore, lemma \ref{covering00} applies with $f=|u|^{p+1}$ (which is $L^{p+1}_{loc}({\cal D}_{u})$  from the solution of the Cauchy problem) and we write from iii) of lemma \ref{covering00}
\begin{multline*}
\sup_{\{x\;|\;|x-x_0|\le
 \frac{e^{-s}}{\delta_0}\}} \left (\int_{s}^{s+10}\!\!
\int_{B} \!\!|w_{x,T^*\!(x)}(y,\tau)|^{p+1} dy d\tau \right)
\\ \le C(\delta_0)
\sup_{\{x\;|\;|x-x_0|\le \frac{e^{-s}}{\delta_0}\}}\left( \int_{s}^{s+10}\!\! \int_{B_{\frac12}} \!\!|w_{x,T^*\!(x)}(y,\tau)|^{p+1}
 dy d\tau\right).
\end{multline*}
Since the right-hand side is bounded by corollary \ref{cor10}, the save holds for the left-hand side (use in particular the Remark 2.2).
This concludes the proof of proposition \ref{tt}.

\Box


\section{Boundedness of the solution in similarity variables}
This section is divided in two parts:
  \begin{itemize}
\item  We first state a general version of  theorem \ref{lyap}, uniform for $x$ near $x_0$ and prove it.
Then, we give a blow-up criterion for  equation (\ref{B}) based on the Lyapunov functional.
    \item  We prove  theorem \ref{t}.
  \end{itemize}
\subsection{A Lyapunov functional }
Consider $u $  a solution of ({\ref{1}}) with blow-up graph $\Gamma:\{x\mapsto T(x)\}$ and  $x_0$ is a non characteristic point.
Consider   $T_0\in
 (t_0(x_0), T(x_0)]$, where $t_0(x_0)$ is defined in (\ref{new1}), for all
$x\in \er^N$ is such that $|x-x_0|\le \frac{T_0-t_0(x_0)}{\delta_0(x_0)}$, where $\delta_0(x_0)$ is defined in (\ref{cone1}), then we write $w$ instead of $w_{x,T^*(x)}$ defined in (\ref{scaling}) with $T^*(x)$ given in  (\ref{explosionx}).
We aim at proving that the functional
  $H$ defined in (\ref{f1}) is a Lyapunov functional  for equation (\ref{B}), provided that $s$ is large enough.
\begin{lem}
\label{lyap0}
For all $ s\ge -\log (T^*(x)-t_0(x_0))$, we have
 the following inequality,
\begin{eqnarray}\label{lem21}
\frac{d}{ds}(E(w))\!\!\!\!&\le&\!\!\!\!-\int_{\partial B}(\partial_s w)^2(\sigma,s){\mathrm{d}}\sigma
+\Sigma (s),
\end{eqnarray}
where $\Sigma(s)$ satisfies
\begin{equation}\label{lem221}
\Sigma (s) \le C e^{-2\gamma s}+ C e^{-2\gamma s}\int_{B}\Big (w^2+|\grad w|^2+(\partial_s w)^2+|w|^{p+1} \Big ){\mathrm{d}}y
 .\qquad \qquad
\end{equation}
\end{lem}
{\it Proof:}  Multipling $(\ref{B})$ by $\partial_s w$, and integrating over the ball $B$, we obtain, for all $s\ge -\log (T^*(x)-t_0(x_0))$,
(recall from \cite{MZ2} that in the case where, $(f,g)\equiv (0,0)$, we have $\frac{d}{ds}E_0(w)=- \int_{\partial B}(\partial_s w)^2(\sigma,s){\mathrm{d}}\sigma$).\\
\begin{eqnarray}\label{E0}
\frac{d}{ds}(E_0(w)+I_0 (w))&=&- \int_{\partial B}(\partial_s w)^2(\sigma,s){\mathrm{d}}\sigma+
\underbrace{\frac{2(p+1)}{p-1}e^{-\frac{2(p+1)s}{p-1}}\int_{B}F\Big(e^{\frac{2s}{p-1}}w\Big) {\mathrm{d}}y}_{I_1} \nonumber\\
&&+\underbrace{\frac{2}{p-1}e^{-\frac{2ps}{p-1}}\int_{B}f\Big(e^{\frac{2s}{p-1}}w\Big)w {\mathrm{d}}y }_{I_2}\nonumber\\
&&+\underbrace{ e^{-\frac{2ps}{p-1}}\int_{B}g\Big(e^{\frac{(p+1)s}{p-1}}(\partial_s w+y.\grad w+\frac{2}{p-1}w)\Big)\partial_s w {\mathrm{d}}y}_{I_3}.
\end{eqnarray}
By exploiting the fact that
$ |{F(x)}|+|x{f(x)}|\le C( 1+|x|^{p+1}),$  we deduce  that
for all $s \ge -\log (T^*(x)-t_0(x_0))$,
\begin{equation}\label{I100}
|I_1|+|I_2|\le C e^{-2\gamma s}+ Ce^{-2\gamma s}\int_{B}|w|^{p+1} {\mathrm{d}}y.
\end{equation}
Note that by combining the inequality
$ |g(x)|\le M( 1+|x|)$ and the fact that  $-\log (T^*(x)-t_0(x_0))\ge 0,$ we obtain
\begin{eqnarray}\label{I122}
|I_3|&\le&
C e^{-2\gamma s}+
 Ce^{-2\gamma s}
\int_{B}\Big ((\partial_s w)^2 +w^2+|\grad w |^2 \Big ) {\mathrm{d}}y.\qquad
\end{eqnarray}
Then, by using  (\ref{E0}), (\ref{I100}) and (\ref{I122}),  we have the estimates (\ref{lem21})  and (\ref{lem221})
for all $s \ge -\log (T^*(x)-t_0(x_0))$. This concludes the proof of lemma
 \ref{lyap0}.
\Box

\vspace{0.1cm}

\no With lemma \ref{lyap0} and  proposition \ref{tt} we are in a position to prove theorem \ref{lyap}'.
\\

\noindent {\bf{THEOREM} 1.1'} {\it (Existence  of a Lyapunov functional for equation }
({\ref{B}}))\\
\label{lyapbis}
{\it  Consider $u $   a solution of ({\ref{1}}) with blow-up graph
 $\Gamma:\{x\mapsto T(x)\}$ and  $x_0$ is a non characteristic point. Then
  there exists  $t_0(x_0)\in [0, T(x_0))$ such that,  for all $T_0\in
 (t_0(x_0), T_0(x_0)]$, for all  $s\ge -\log (T_0-t_0(x_0))$
and $x\in \er^N$, where $|x-x_0|\le \frac{e^{-s}}{\delta_0(x_0)}$,
we have
\begin{eqnarray}\label{L1}
H(w(s+10))-H(w(s))
&\le&  - \int_{s}^{s+10}\int_{\partial B}(\partial_sw)^2(\sigma,\tau)d\sigma d\tau, \qquad \qquad
\end{eqnarray}
where $w=w_{x,T^*(x)}$ and  $T^*(x)$ is defined in (\ref{explosionx}).}


\bigskip

\no  {\it{Proof of theorem \ref{lyap}':}}
Consider $u$  is a solution of ({\ref{1}}) with blow-up graph $\Gamma:\{x\mapsto T(x)\}$ and  $x_0$ is a non characteristic point. We can apply  proposition \ref{tt}, if $\m=\m_0$ small enough, where we have $\frac{\m(p+3)}{2}\le \gamma$, 
for all $T_0\in
 (t_0(x_0), T(x_0)]$, for all $s\ge -\log (T_0-t_0(x_0))$ and
$x\in \er^N$, where $|x-x_0|\le \frac{e^{-s}}{\delta_0(x_0)}$,  we get
\begin{eqnarray}\label{Az}
\int_{s}^{s+10} \!\int_{B}\!(\partial_s w)^2(y,\tau){\mathrm{d}}y{\mathrm{d}}\tau
+\int_{s}^{s+10} \int_{B}|w(y,\tau)|^{p+1}{\mathrm{d}}y{\mathrm{d}}\tau\nonumber\\
 + \int_{s}^{s+10} \int_{B}|\grad w(y,\tau)|^2{\mathrm{d}}y{\mathrm{d}}\tau \le K_1e^{\gamma s}.
\end{eqnarray}
By combining (\ref{Az}) and   the  following Jensen's inequality
\begin{equation*}
\int_{B}w^2 {\mathrm{d}}y\le C
 \int_{B}|w|^{p+1}{\mathrm{d}}y+C,
\end{equation*}
we obtain
\begin{eqnarray}\label{I1222}
 \int_{s}^{s+10}e^{-2\gamma \tau}
\int_{B}\Big ((\partial_s w)^2 +w^2+|\grad w |^2+|w|^{p+1} \Big ) {\mathrm{d}}y{\mathrm{d}}\tau
&\le&
 CK_1e^{-\gamma s}.
\end{eqnarray}
Therefore, lemma \ref{lyap0}
 and  (\ref{I1222}), implies
\begin{eqnarray}\label{}
E(w)(s+10)-E(w)(s)\!\!\!\!&\le&\!\!\!\!-\int_{s}^{s+10}\int_{\partial B}(\partial_s w)^2(\sigma,\tau){\mathrm{d}}\sigma{\mathrm{d}}\tau
+ CK_1e^{-\gamma s}.
\end{eqnarray}
Then, we write
\begin{eqnarray*}\label{H002}
H(w)(s+10)-H(w)(s)\le -\int_{s}^{s+10}\!\int_{\partial B}\!(\partial_s w)^2(\sigma,\tau){\mathrm{d}}\sigma{\mathrm{d}}\tau
+ (CK_1-\sigma (1- e^{-10\gamma }))
 e^{-\gamma s}, \qquad
\end{eqnarray*}
where $H$ is defined in (\ref{f1}).

\no We now choose $\sigma$ large enough,  so we have $CK_1-\sigma (1- e^{-10\gamma })\le 0$ and then
 (\ref{L1}) is a direct consequence of the above inequality. 
This concludes the proof of theorem \ref{lyap}'.

\Box

\no We now claim the following lemma:
\begin{lem}\label{L02}
 There exists $S_2\ge S_1$ such that, if
  $H(w(s_2))<0$ for some $s_2\ge \max(S_2,-\log (T^*(x)-t_0(x_0)))$, then $w$
blows up in some finite time $S>s_2$.
\end{lem}
\no {\it Proof:}
The argument is the same as in the corresponding part
in \cite{AM}.

\Box

\vspace{0.1cm}

\subsection{Boundedness of the solution in similarity variables}
\no We prove theorem 1.2 here. Note that the lower bound follows from the finite speed of propagation and wellposedness  in $H^1\times L^2$.
For a detailed argument in the similar case of equation (\ref{mu=0}), see lemma 3.1 (page 1136) in \cite{MZ3}.\\

\no We define the following time:
\begin{equation}\label{new2}
 t_1(x_0)=\max(T(x_0)-e^{-S_2},0).
\end{equation}
Given some $T_0\in (t_1(x_0),T(x_0)]$, for all
$x\in \er^N$ is such that $|x-x_0|\le \frac{T_0-t_1(x_0)}{\delta_0(x_0)}$, where $\delta(x_0)$ is defined in (\ref{cone1}), then we write $w$ instead of $w_{x,T^*(x)}$ defined in (\ref{scaling}) with $T^*(x)$ given in  (\ref{explosionx}).
We aim   at bounding $\|(w,\partial_s w)(s)\|_{H^1\times L^2(B)}$ for $s$ large.

\no As in \cite{MZ2} , by combining theorem \ref{lyap} and  lemma \ref{L02} we get the following bounds:
\begin{cor} (Bounds on $E$) For all $s\ge -\log (T^*(x)-t_1(x_0))$,  it holds that
\label{cor1}
\begin{eqnarray}\label{cor01}
-C\le E(w(s))\le  M_0\nonumber\\
\int^{s+10}_{s}\int_{\partial B}\Big (\partial_s w(\sigma,\tau)\Big )^2d\sigma d\tau\le M_0,
\end{eqnarray}
where  $M_0=M_0(T_0,\|(u(t_0(x_0)),\partial_tu(t_0(x_0)))\|_{
H^{1}\times L^{2}(B(x_0,\frac{T(x_0)-t_0(x_0)}{\delta_0(x_0)})
)})$, $C>0$
   and   $\delta_0(x_0)\in (0,1)$ is defined in (\ref{cone1}).
\end{cor}
We now claim the following corollary
\begin{cor}  For all  $s\ge  -\log (T^*(x)-t_1(x_0))$, we have
\label{cor100}
\begin{eqnarray}\label{cor00001}
\int_{s}^{s+10}\int_{B}\Big (\partial_s w(y,\tau)-\lambda(\tau,s) w(y,\tau)\Big )^2dy d\tau\le CM_0,
\end{eqnarray}
and (nonconcentration property) for all $b\in \R^N$ and $r_0\in (0,1)$ such that
$B(b,r_0)\subset B(0,\frac{1}{\delta_0(x_0)})$,
\begin{eqnarray}\label{cor00001bis}
\int_{s}^{s+\sqrt{r_0}}\int_{B_{(b,r_0)}}\Big (\partial_s w(y,\tau)-\lambda(\tau,s) w(y,\tau)\Big )^2dy d\tau\le CM_0 r_0,
\end{eqnarray}
where $0\le \lambda(\tau,s) \le C(\delta_0)$, for all $\tau \in [s,s+\sqrt{r_0}]$.
\end{cor}
\no {\it Proof:}
The argument is the same as in the corresponding part, see proposition 4.2 (page 1147)
in   \cite{MZ3}.

\Box

\vspace{0.1cm}

\no Starting from these bounds, the proof of theorem 1.2 is similar to the proof in \cite{MZ3}
except for
the treatment
of the perturbation terms. In our opinion, handling these terms is straightforward in
all the steps of the proof, except for the first step, where we bound the time averages of the $L^{p+1}(B)$ norm of $w$. For that reason,
we only give that step and refer to  \cite{MZ3} for the remaining steps in the proof of  theorem 1.2. This is the step we prove here (In the following
$K_2$, $K_3$ denotes  a constant that depends only on $C$, $M_0$, $\varepsilon_1$ and $\varepsilon_2$
is an arbitrary positive number in $(0,1)$).
\begin{pro}\label{pro}(Control of the space-time $L^{p+1}$ norm of $w$)\\
 For all  $s\ge -\log (T^*(x)-t_2(x_0))$,  for some $t_2(x_0)\in [t_1(x_0),T(x_0))$, for all $\varepsilon_1>0$,
 \begin{eqnarray}\label{pro1}
\int_{s}^{s+10}\!\!\!\int_{B}\!| w|^{p+1}{\mathrm{d}}y{\mathrm{d}}\tau&\le& \frac{K_3}{\varepsilon_1}+
+K_3\varepsilon_1
\int_{s}^{s+10}\!\!\!\int_{B}\!|\grad w|^2{\mathrm{d}}y{\mathrm{d}}\tau\nonumber\\
&&+
C \int_{B}\!\Big(| \partial_sw(y,s)|^2+| \partial_sw(y,s+10)|^2\Big){\mathrm{d}}y
\end{eqnarray}
\end{pro}
{\it Proof:} By integrating the expression (\ref{f2}) of $E$ in time between $s$ and $s+10$, where $s\ge -\log (T^*(x)-t_1(x_0))$, we obtain:
\begin{eqnarray}\nonumber
\int_{s}^{s+10}\!\!E(\tau)  d\tau&=&\displaystyle\int_{s}^{s+10}\int_{B}\!\!\Big ( \frac{1}{2}(\partial_s w)^2
+\frac{p+1}{(p-1)^2}w^2-\frac{1}{p+1}|w|^{p+1}\Big ) {\mathrm{d}}y{\mathrm{d}}\tau\\
&&+\frac{1}{2}\displaystyle\int_{s}^{s+10}\!\!\int_{B}\!\!\Big (|\grad w|^2-(y.\grad w)^2\Big ) {\mathrm{d}}y{\mathrm{d}}\tau \nonumber\\
&&
-\int_{s}^{s+10}\!\! e^{-\frac{2(p+1)\tau}{p-1}}\displaystyle\int_{B}F(e^{\frac{2}{p-1}\tau}w) {\mathrm{d}}y\mathrm{d}\tau.\label{et}
\end{eqnarray}
By multiplying the equation (\ref{B}) by $w$ and integrating both in time and in space over $B\times [s,s+10]$, we obtain
the following identity, after some
integration by parts:
\begin{eqnarray}\label{et1}
&&\Big [\int_{B}\!\!\Big (w\partial_s w+(\frac{p+3}{2(p-1)}-N)w^2\Big ) {\mathrm{d}}y\Big ]_{s}^{s+10}=
\int_{s}^{s+10}\!\!\int_{B}\!\!(\partial_s w)^2{\mathrm{d}}y{\mathrm{d}}\tau\nonumber\\
&&-\int_{s}^{s+10}\!\!\int_{B}\!\!(|\grad w|^2-(y.\grad w)^2){\mathrm{d}}y{\mathrm{d}}\tau
-\frac{2p+2}{(p-1)^2}\int_{s}^{s+10}\!\!\int_{B}\!\!w^2{\mathrm{d}}y{\mathrm{d}}\tau\nonumber\\
&&+\int_{s}^{s+10}\!\!\int_{B}\!\!|w|^{p+1}{\mathrm{d}}y{\mathrm{d}}\tau
+2\!\!\int_{s}^{s+10}\!\!\!\int_{B}\!\!\!\partial_s w(y.\grad w) {\mathrm{d}}y{\mathrm{d}}\tau-\!2\!\int_{s}^{s+10}\!\!\int_{\partial B}
\!\!w\partial_s w {\mathrm{d}}\sigma{\mathrm{d}}\tau\nonumber\\
&&+\!\int_{s}^{s+10}\!\!\int_{B}\!\!e^{-\frac{2p\tau}{p-1}}f\Big(e^{\frac{2\tau}{p-1}}w\Big)w{\mathrm{d}}y{\mathrm{d}}\tau\nonumber\\
&& +\!\int_{s}^{s+10}\!\!\int_{B}\!\!
e^{-\frac{2p\tau}{p-1}}g\Big(e^{\frac{(p+1)\tau}{p-1}}(\partial_s w+y.\grad w+\frac{2}{p-1}w)\Big)w{\mathrm{d}}y{\mathrm{d}}\tau.\qquad
\end{eqnarray}
By combining the identities (\ref{et}) and (\ref{et1}), we obtain
\begin{eqnarray}\label{controlp}
\frac{(p-1)}{2(p+1)}&&\int_{s}^{s+10}\!\!\int_{B}\!\!|w|^{p+1}{\mathrm{d}}y{\mathrm{d}}\tau = \nonumber\\
&&\frac{1}{2}\Big [\int_{B}\!\!\Big (w\partial_s w+(\frac{p+3}{2(p-1)}-N)w^2\Big ) {\mathrm{d}}y\Big ]_{s}^{s+10}-
\int_{s}^{s+10}\!\!\int_{B}\!\!(\partial_s w)^2{\mathrm{d}}y{\mathrm{d}}\tau\nonumber\\
&&+\int_{s}^{s+10}\!\!E(s) d\tau
-\int_{s}^{s+10}\!\!\int_{B}\!\!\partial_s w(y.\grad w) {\mathrm{d}}y{\mathrm{d}}\tau+\!\int_{s}^{s+10}\!\!\int_{\partial B}\!\!w\partial_s w {\mathrm{d}}\sigma{\mathrm{d}}\tau\nonumber\\
&& -\underbrace{\frac12\!\int_{s}^{s+10}\!\!\int_{B}\!\!
e^{-\frac{2p\tau}{p-1}}g\Big(e^{\frac{(p+1)\tau}{p-1}}(\partial_s w+y.\grad w+\frac{2}{p-1}w)\Big)w{\mathrm{d}}y
{\mathrm{d}}\tau}_{A_1}\nonumber\\
&&-\underbrace{\frac12\!\int_{s}^{s+10}\!\!\int_{B}\!\!e^{-\frac{2p\tau}{p-1}}f\Big(e^{\frac{2\tau}{p-1}}w
\Big)w{\mathrm{d}}y{\mathrm{d}}\tau}_{A_2}\nonumber\\
&&+\underbrace{\int_{s}^{s+10}\!\! e^{-\frac{2(p+1)\tau}{p-1}}\displaystyle\int_{B}F(e^{\frac{2}{p-1}\tau}w)
{\mathrm{d}}y\mathrm{d}\tau}_{A3}.
\end{eqnarray}
We claim that proposition \ref{pro} follows from the following lemma where we   have the following estimates
\begin{lem}\label{g}
 For all  $ s\ge -\log (T^*(x)-t_2(x_0))$,  for some $t_2(x_0)\in [t_1(x_0),T(x_0))$,
 for all $\varepsilon_1>0$, for all $\varepsilon_2>0$,
\begin{equation}\label{control}
\int_{s}^{s+10}\int_{B}\!\!|\grad w|^2(1-|y|^2){\mathrm{d}}y{\mathrm{d}}\tau\le K_2 +C
\int_{s}^{s+10}\!\!\int_{B}\!\!| w|^{p+1}{\mathrm{d}}y{\mathrm{d}}\tau,
\end{equation}
\begin{equation}\label{control1}
\sup_{\tau \in [s,s+10]}\int_{B} \!\!w^2(y,\tau){\mathrm{d}}y\le \frac{K_2}{\varepsilon_2} +K_2\varepsilon_2
\int_{s}^{s+10}\!\!\int_{B}\!\!|w|^{p+1}{\mathrm{d}}y{\mathrm{d}}\tau,
\end{equation}
\begin{eqnarray}\label{control3}
\int_{s}^{s+10}\!\!\!\int_{B}\!|\partial_s w y.\grad w| {\mathrm{d}}y{\mathrm{d}}\tau
&\le& \frac{K_2}{\varepsilon_1\varepsilon_2} +K_2\varepsilon_1
\int_{s}^{s+10}\!\!\!\int_{B}\!|\grad w|^2{\mathrm{d}}y{\mathrm{d}}\tau\nonumber\\
&& +K_2\varepsilon_2
\int_{s}^{s+10}\!\!\int_{B}\!\!|w|^{p+1}{\mathrm{d}}y{\mathrm{d}}\tau,
\end{eqnarray}
\begin{eqnarray}\label{control30}
\int_{s}^{s+10}\!\!\!\int_{\partial B}\!|\partial_s ww  |{\mathrm{d}}\sigma{\mathrm{d}}\tau
&\le& \frac{K_2}{\varepsilon_1 \varepsilon_2} +K_2\varepsilon_2
\int_{s}^{s+10}\!\!\!\int_{B}\!| w|^{p+1}{\mathrm{d}}y{\mathrm{d}}\tau\nonumber\\
&&+K_2\varepsilon_1
\int_{s}^{s+10}\!\!\!\int_{B}\!|\grad w|^2{\mathrm{d}}y{\mathrm{d}}\tau,\qquad\qquad\qquad
\end{eqnarray}
\begin{equation}\label{control4}
\int_{B}|w\partial_s w|{\mathrm{d}}y\le
\int_{B} (\partial_s w)^2{\mathrm{d}}y+\frac{K_2}{\varepsilon_2} +K_2\varepsilon_2
\int_{s}^{s+10}\!\!\!\int_{B}\!| w|^{p+1}{\mathrm{d}}y{\mathrm{d}}\tau,
\end{equation}
\begin{eqnarray}\label{A10}
|A_1|&\le& \frac{K_2}{\varepsilon_2} +K_2\varepsilon_2
\int_{s}^{s+10}\!\!\!\int_{B}\!| w|^{p+1}{\mathrm{d}}y{\mathrm{d}}\tau+
C  e^{-2s}\int_{s}^{s+10}\!\!\!\int_{B}\!|\grad w|^2{\mathrm{d}}y{\mathrm{d}}\tau,
\end{eqnarray}
\begin{eqnarray}\label{A20}
|A_2|+|A_3|&\le&  C+ Ce^{-2\gamma s}\int_{s}^{s+10} \!\!\int_{B}|w|^{p+1} {\mathrm{d}}y{\mathrm{d}}\tau.
\end{eqnarray}
\end{lem}
Indeed, from (\ref{controlp}) and this lemma,
we deduce that
\begin{eqnarray*}
\int_{s}^{s+10}\!\!\!\int_{B}\!| w|^{p+1}{\mathrm{d}}y{\mathrm{d}}\tau&\le& \frac{K_2}{\varepsilon_1\varepsilon_2}+
(K_2\varepsilon_2+C  e^{-2\gamma s})
\int_{s}^{s+10}\!\!\!\int_{B}\!| w|^{p+1}{\mathrm{d}}y{\mathrm{d}}\tau\\
&&+(K_2\varepsilon_1+
C  e^{-2s})
\int_{s}^{s+10}\!\!\!\int_{B}\!|\grad w|^2{\mathrm{d}}y{\mathrm{d}}\tau\\
&&+
C \int_{B}\!\Big(| \partial_sw(y,s)|^2+| \partial_sw(y,s+10)|^2\Big){\mathrm{d}}y
\end{eqnarray*}
Now, we can use the fact that
$s\ge -\log(T^*(x)-t_2(x_0))\ge -\log(T(x_0)-t_2(x_0))$ and we prove easily that
 $$e^{-2s}\le (T(x_0)-t_2(x_0))^2\qquad and\qquad
 e^{-2\gamma s}\le (T(x_0)-t_2(x_0))^{2\gamma}. $$
Taking $t_2(x_0)$, where $ T(x_0)-t_2(x_0)$ small enough, where we have
$$ C(T(x_0)-t_2(x_0))^2\le K_2\varepsilon_1\qquad and\qquad
  (T(x_0)-t_2(x_0))^{2\gamma}\le K_2\varepsilon_2. $$
 If we choose $\varepsilon_2=\frac{1}{4K_2}$, we obtain (\ref{pro1}).

\bigskip

\no
It remains to prove lemma 3.6.\\

\bigskip

\no
{\it{ Proof of lemma 3.6:}}
 For the estimates (\ref{control}), (\ref{control1}), (\ref{control3}), (\ref{control30})  and (\ref{control4}),  we can adapt with no difficulty
 the proof given in the case of  the wave equation treated in \cite{MZ1}.

\no
Now, we control the terms $A_1$, $A_2$ and $A_3$. Since
$ |g(x)|\le M(1+|x|)$, we have
\begin{eqnarray}
\label{I12}
|A_1|&\le& Ce^{-s}\int_{s}^{s+10}\!\!
\int_{B}(\partial_s w)^2 {\mathrm{d}}y{\mathrm{d}}\tau +
Ce^{-s}\int_{s}^{s+10}\!\!\int_{B}|y.\grad w ||w| {\mathrm{d}}y{\mathrm{d}}\tau\nonumber\\
&&+ Ce^{-s} \int_{s}^{s+10}\!\!\int_{B}w^2 {\mathrm{d}}y{\mathrm{d}}\tau+ C e^{-\frac{2p s}{p-1}}\int_{s}^{s+10}\!\!\int_{B}|w| {\mathrm{d}}y{\mathrm{d}}\tau.\qquad
\end{eqnarray}
By using (\ref{cor00001}), we get
\begin{eqnarray}\label{CC0}
 \int_{s}^{s+10}\!\!
\int_{B}(\partial_s w)^2 {\mathrm{d}}y{\mathrm{d}}\tau \le K_2+C\int_{s}^{s+10}\!\!
\int_{B}w^2 {\mathrm{d}}y{\mathrm{d}}\tau.
\end{eqnarray}
As usual, we write
\begin{eqnarray}
 \label{CC2}
C  \int_{s}^{s+10}\int_{B}|w| {\mathrm{d}}y{\mathrm{d}}\tau&\le& C  +C\int_{s}^{s+10}  \int_{B}w^2 {\mathrm{d}}y{\mathrm{d}}\tau.
\end{eqnarray}
Using the Cauchy-Schwarz inequality, we obtain
\begin{eqnarray}\label{controlp3}
C\!e^{-s} \int_{s}^{s+10}\!\!\!\!\int_{B}\!\!|w||y.\grad w|{}{\mathrm{d}}y{\mathrm{d}}\tau
\!\!&\le&\!\!   \int_{s}^{s+10}\!\!\int_{B}\!\!w^{2}{\mathrm{d}}y{\mathrm{d}}\tau
+C e^{-2s}\!\! \int_{s}^{s+10}\!\!\int_{B}\!\!|\grad w|^2{\mathrm{d}}y{\mathrm{d}}\tau.\qquad
\end{eqnarray}
By combining the fact that $ s\ge0$ and the  inequalities (\ref{control1}),
(\ref{I12}), (\ref{CC0}),  (\ref{CC2})  (\ref{controlp3}), we deduce
\begin{eqnarray}\label{A1}
|A_1|&\le& \frac{K_2}{\varepsilon_2} +K_2\varepsilon_2\int_{s}^{s+10}\!\!\!\int_{B}\!| w|^{p+1}{\mathrm{d}}y{\mathrm{d}}\tau
+C  e^{-2s}\int_{s}^{s+10}\!\!\!\int_{B}\!|\grad w|^2{\mathrm{d}}y{\mathrm{d}}\tau.
\end{eqnarray}
Similarly,
by exploiting the fact that
$ |{F(x)}|+   |xf(x)|\le C(1+|x|^{p+1}),$
we obtain
\begin{eqnarray}\label{A3}
|A_2|+|A_3|&\le&  C\int_{s}^{s+10} \!\!e^{-2\gamma \tau}{\mathrm{d}}\tau+ C\int_{s}^{s+10} \!\!e^{-2\gamma \tau}\int_{B}|w|^{p+1} {\mathrm{d}}y{\mathrm{d}}\tau\nonumber\\
&\le&  C+ Ce^{-2\gamma s}\int_{s}^{s+10} \!\!\int_{B}|w|^{p+1} {\mathrm{d}}y{\mathrm{d}}\tau.
\end{eqnarray}
This concludes the proof of lemma 3.6 and proposition 3.5 too.
\Box

\no Since the derivation of theorem 1.2  from proposition 3.5  is the
same as in the non perturbed case treated in \cite{MZ3} (up to some very minor changes),
this concludes the proof of theorem 1.2.
\Box

\appendix

\section{Covering technique}
 In this section, we prove the covering result stated in lemma \ref{covering00}.


\no {\it Proof of lemma \ref{covering00}}:\\
(i) For all $x$ such that $|x-x_0|\le \frac{T_0-t_1}{\delta_0}$, we easily have
\begin{equation}\label{merle}
T^*(x)-t_1\le T_0-t_1,\  t_2(x)\le t_2(x_0) \mbox{ and }B(x, T^*(x)-t_1)\subset
B\left(x_0,  \frac{T_0-t_1}{\delta_0}\right),
\end{equation}
the basis of ${\cal D}_{x_0, T_0, t_1, \delta_0}$.
Therefore
\begin{equation}\label{ap1}
{\cal S}_{x, T^*(x),t_1,1}\subset
{\cal S}_{x_0, T_0,t_1, \delta_0} \subset {\cal D}_u
 \end{equation}
and  (i)  follows.
\medskip

\noindent (ii) Since $T^*(x)\le T_0$ and $t_2(x)\le t_2(x_0)$, we write from (\ref{ap1})
\[
\int_{ {\cal S}_{x,T^*(x),t_1,1}} (T^*(x)-t)^\kappa |f(\xi,\t0)|^q d\xi dt\le
\int_{ {\cal S}_{x_0,T_0,t_1,\delta_0}} (T_0-t)^\kappa |f(\xi,\t0)|^q d\xi dt
\]
which is finite since $f\in L^q_{loc}\left(D_u\right)$ and because there exists $0<r_0<1$ such that ${\cal S}_{x_0,T_0,t_1,r_0\delta_0}\subset D_u$ (this is true because in (\ref{ap1}), $D_u$ is open and ${\cal S}_{x_0,T_0,t_1,\delta_0}$ is closed). Thus, the supremum exists and (ii) is proved.

\medskip

\noindent (iii) Consider $x^*$ such that $|x^*-x_0|\le \frac{T_0-\t0}{\delta_0}$ and
\[\int_{ {\cal S}_{x^*,T^*(x^*),t_1,1}}\!(T^*(x^*)-\t0)^\kappa  |f(\xi,\t0)|^q d\xi dt\ge \frac 12 \sup_{|x-x_0|\le \frac{T_0-\t0}{\delta_0}}
\int_{ {\cal S}_{x,T^*(x),t_1,1}}
(T^*(x)-\t0)^\kappa |f(\xi,\t0)|^q d\xi dt.
\]
It is enough to prove that
\begin{eqnarray}\label{supineq}
\int_{ {\cal S}_{x^*,T^*(x^*),t_1,1}}
(T^*(x^*)-t)^\kappa  |f(\xi,\t0)|^q d\xi dt
 \qquad\qquad\qquad\qquad\qquad\qquad\qquad\\
\le C(\delta_0,\kappa) \sup_{|x-x_0|\le \frac{T_0-t_1
}{\delta_0}}\int_{t_1}^{t_2(x)} (T^*(x)-\t0)^\kappa \int_{B(x, \frac{T^*(x)-\t0}2)} |f(\xi,\t0)|^q d\xi dt\nonumber
\end{eqnarray}
in order to conclude. In the following, we will prove \aref{supineq}.

\medskip

Note that we can cover
${\cal S}_{x^*,T^*(x^*),t_1,1} $
by $k(\delta_0)$ sets of the form
${\cal S}_{x_i,\widetilde{T^*}(x_i),t_1,\frac{1-\delta_0}2}$
 where $|x_i-x^*|\le T^*(x^*)-t_1$ and
 $$\widetilde{T^*}(x_i)
={T^*}(x^*)- \delta_0|x_i-x^*|$$
is the later boundary in the backward cone of vertex $(x^*,T^*(x^*))$ and slape  $\delta_0$.
Indeed, this number does not change by scaling and is thus the same as the number when ${T^*}(x^*)-t_1=1$.
In addition,
\[
\begin{array}{cl}
&\left|(T^*(x_i)-t_1)-(T^*(x^*)-t_1)\right|= \left|T^*(x_i)-T^*(x^*)\right|
=\delta_0\left||x_i-x_0|-|x^*-x_0|\right|\\
\\
\le & \delta_0|x_i-x^*|\le \delta_0 (T^*(x^*)-t_1),
\end{array}
\]
hence
\begin{equation}\label{app3}
(1-\delta_0)(T^*(x^*)-t_1)\le T^*(x_i)-t_1\le (1+\delta_0)(T^*(x^*)-t_1).
\end{equation}
Since $\widetilde{T^*}(x_i)\le T^*(x_i)$, it follows that
$${\cal S}_{x^*,T^*(x^*),t_1,1} \subset \bigcup_{i=1}^{k(\delta_0)}
{\cal S}_{x_i,\widetilde{T^*}(x_i),t_1,\frac{1-\delta_0}2}
\subset \bigcup_{i=1}^{k(\delta_0)}\Big(
{\cal S}_{x_i,T^*(x_i),t_1,\frac12}\bigcap \{(\xi,\tau )|\ \tau\le t_2(x^*)\}\Big).$$
Moreover, for any $i=1,...k(\delta_0)$ and $t\in [t_1, \min (t_2(x^*),t_2(x_i))]$, we write
$$e^{-10}
( T^*(x_i)-t_1)= T^*(x_i)-t_2(x_i)\le
 T^*(x_i)-t\le T^*(x_i)-t_1,
$$
$$e^{-10}
( T^*(x^*)-t_1)= T^*(x^*)-t_2(x^*)\le
 T^*(x^*)-t\le T^*(x^*)-t_1.
$$
Using (\ref{app3}), we see that
$$e^{-10}(1-\delta_0)
( T^*(x^*)-t)\le
 T^*(x_i)-t\le e^{10}(1+\delta_0)(T^*(x^*)-t).
$$
It follows then that
\begin{eqnarray*}
&&
\int_{
{\cal S}_{x^*,T^*(x^*),t_1,1}}
(T^*(x^*)-\t0)^\kappa  |f(\xi,\t0)|^q d\xi dt\\
&\le&  \sum_{i=1}^{k(\delta_0)} \int_{
{\cal S}_{x_i,T^*(x_i),t_1,\frac12}
\cap \{(\xi,t ),\ t\le t_2(x^*)\}
} (T^*(x^*)-\t0)^\kappa |f(\xi,\t0)|^q d\xi dt\\
&\le & \sum_{i=1}^{k(\delta_0)}
\int_{
{\cal S}_{x_i,T^*(x_i),t_1,\frac12}}
 \frac{e^{10\kappa}}{(1-\delta_0)^\kappa}(T^*(x_i)-\t0)^\kappa |f(\xi,\t0)|^q d\xi dt\\
&\le & \frac{k(\delta_0)e^{10\kappa}}{(1-\delta_0)^\kappa} \sup_{|x-x_0|\le \frac{T_0-t_1}{\delta_0}} \int_{
{\cal S}_{x,T^*(x),t_1,\frac12}} (T^*(x)-t)^\kappa  |f(\xi,\t0)|^q d\xi dt,
\end{eqnarray*}
where we used in the last line the fact that $x_i\in B(x^*, T^*(x^*)-t_1)\subset B(x_0, \frac{T_0-t_1}{\delta_0})$ by \aref{merle}.
This yields (iii) and concludes the proof of lemma \ref{covering00}.\Box

\noindent{\bf Address}:\\
Universit\'e de Tunis El-Manar, Facult\'e des Sciences de Tunis, D\'epartement de math\'ematiques, Campus Universitaire 1060,
 Tunis, Tunisia.\\
\vspace{-7mm}
\begin{verbatim}
e-mail: ma.hamza@fst.rnu.tn
\end{verbatim}
Universit\'e Paris 13, Institut Galil\'ee,
Laboratoire Analyse, G\'eom\'etrie et Applications, CNRS UMR 7539,
99 avenue J.B. Cl\'ement, 93430 Villetaneuse, France.\\
\vspace{-7mm}
\begin{verbatim}
e-mail: Hatem.Zaag@univ-paris13.fr
\end{verbatim}

\end{document}